\documentclass{amsart}
\usepackage{amssymb}

\usepackage{amsmath}
\usepackage{amscd}
\usepackage{thmdefs}

\input{tcilatex}
\theoremstyle{definition}
\theoremstyle{remark}
\numberwithin{equation}{section}

\input tcilatex

\begin{document}
\title[Subalgebras of $\left( W^{*}(G),E\right) $]{Subalgebras and Free Product Structures of a Graph $W^{*}$-Probability Space}
\author{Ilwoo Cho}
\address{Univ. of Iowa, Dep. of Math, Iowa City, IA, U. S. A}
\email{ilcho@math.uiowa.edu}
\date{}
\dedicatory{}
\thanks{}
\maketitle

\begin{abstract}
Let $G$ be a countable directed graph. Then we can construct the graph $%
W^{*} $-algebra $W^{*}(G)$ and its diagonal subalgebra $D_{G}.$ By defining
the conditional expectation $E:W^{*}(G)\rightarrow D_{G},$ we have the
grarph $W^{*}$-probability space over $D_{G},$ $\left( W^{*}(G),E\right) ,$
as an amalgamated $W^{*}$-probability space over $D_{G}.$ The amalgamated
freeness on $\left( W^{*}(G),E\right) $ is defined in the sense of Speicher.
In this papaer, we will define the $D_{G}$-semicircular system and $D_{G}$%
-valued R-diagonal system. If the graph $G$ contains $N$-mutually \textbf{%
diagram-distinct} loops, then we can construct the $D_{G}$-semicircular
system $L_{N}$ in $(W^{*}(G),E)$ and the $D_{G}$-semicircular algebra $%
W^{*}(L_{N},D_{G})$ generated by $L_{N}$ and $D_{G}$, characterized by $%
\left( W^{*}(L_{N},D_{N}),E_{N}\right) \otimes \left( D_{G},\mathbf{1}%
\right) ,$ where $D_{N}\leq D_{G}$ is a subalgebra determined by $L_{N}.$ If
we have $N$-mutually distinct finite paths $w_{1},...,w_{N}\notin loop(G),$
then we can construct the $D_{G}$-valued R-diagonal system $%
R=\{L_{w_{1}},L_{w_{1}}^{*},...,L_{w_{N}},L_{w_{N}}^{*}\}$ satisfying that $%
\{L_{w_{1}},L_{w_{1}}^{*}\},...,\{L_{w_{N}},L_{w_{N}}^{*}\}$ are $D_{G}$%
-free from each other. The subalgebra $W^{*}(R,D_{G})$ is observed.
Precisely, we show that the graph $W^{*}$-algebra $W^{*}(G)$ is the $D_{G}$%
-free product of $D_{G}$-free building blocks $D_{G}$ and $%
W^{*}(\{L_{e}\},D_{G}),$ for all $e\in E(G).$ So, we can observe the $D_{G}$%
-free product structure of $\left( W^{*}(G),E\right) $ by these building
blocks. Also, we can see that

$\left( W^{*}(G),E\right) =\left( D_{G},\mathbf{1}\right) *_{D_{G}}(%
\underset{e\in E(G)}{\,*_{D_{G}}}\left( W^{*}(\{L_{e}\},D_{G}),E\right) )$

$\ $
\end{abstract}

\strut

In [16], we constructed the graph $W^{*}$-probability spaces. The graph $%
W^{*}$-probability theory is one of the good example of Speicher's
combinatorial free probability theory with amalgamation. In [16], we
observed how to compute the certain operator-valued moments and cumulants of
an arbitrary operator-valued random variables in the graph $W^{*}$%
-probability space and observed the amalgamated\ freeness on the graph $%
W^{*} $-probability space, with respect to the given conditional
expectation. Also, in [17], we consider certain operator-valued random
variables of the graph $W^{*}$-probability space, for example, semicircular
elements, even elements and R-diagonal elements. This shows that the graph $%
W^{*}$-probability spaces contain the rich free probabilistic objects.

\strut

Throughout this paper, let $G$ be a countable directed graph and let $%
\mathbb{F}^{+}(G)$ be the free semigroupoid of $G,$ in the sense of Kribs
and Power. i.e., it is a collection of all vertices of the graph $G$ as
units and all admissible finite paths, under the admissibility. The
admissible product between two elements in the set $\Bbb{F}^{+}(G)$ is the
binary operation on $\Bbb{F}^{+}(G).$ As a set, the free semigroupoid $%
\mathbb{F}^{+}(G)$ can be decomposed by

\strut

\begin{center}
$\mathbb{F}^{+}(G)=V(G)\cup FP(G),$
\end{center}

\strut

where $V(G)$ is the vertex set of the graph $G$ and $FP(G)$ is the set of
all admissible finite paths. Trivially the edge set $E(G)$ of the graph $G$
is properly contained in $FP(G),$ since all edges of the graph can be
regarded as finite paths with their length $1.$ Kribs and Power defined the
graph Hilbert space $H_{G}$ $=$ $l^{2}$ $\left( \Bbb{F}^{+}(G)\right) $ with
its Hilbert basis $\{\xi _{w}$ $:$ $w$ $\in $ $\Bbb{F}^{+}(G)\}.$ In [16]
and [17], we defined the creation operator $L_{w}$ by the multiplication
operator with its symbol $\xi _{w}$ and its adjoint $L_{w}^{*},$ the
annihilation operator on $H_{G}.$ They have the following relation; if $w$ $%
= $ $v_{1}$ $w$ $v_{2}$ in $\Bbb{F}^{+}(G)$ with $v_{1},$ $v_{2}$ $\in $ $%
V(G), $ then

\strut

(i) \ \ \ $\ L_{w}=L_{v_{1}}L_{w}L_{v_{2}}$

(ii) \ \ $\ L_{w}L_{w}^{*}=L_{v_{1}}$ \ \ and \ \ $L_{w}^{*}L_{w}=L_{v_{2}}$

(iii) \ $\ L_{w}^{2}=L_{w}=L_{w}^{*},$ \ if $w=www\in V(G)$

(iv) \ $\ L_{w}L_{w}^{*}L_{w}=L_{w}$, \ if $w\in FP(G).$\strut

\strut \strut

We define a graph $W^{*}$-algebra of $G$ by

\strut

\begin{center}
$W^{*}(G)\overset{def}{=}\overline{%
\mathbb{C}[\{L_{w},L_{w}^{*}:w\in
\mathbb{F}^{+}(G)\}]}^{w}.$
\end{center}

\strut

Notice that the creation operators induced by vertices are projections and
the creation operators induced by finite paths are partial isometries. We
can define the $W^{*}$-subalgebra $D_{G}$ of $W^{*}(G),$ which is called the
diagonal subalgebra by

\strut

\begin{center}
$D_{G}\overset{def}{=}\overline{\mathbb{C}[\{L_{v}:v\in V(G)\}]}^{w}.$
\end{center}

\strut

Then each element $a$ in the graph $W^{*}$-algebra $W^{*}(G)$ is expressed by

\strut

(0.1) $\ \ \ \ \ \ \ a=\underset{v\in V(G)}{\sum }p_{v}L_{v}+\underset{w\in
FP(G)}{\sum }\left( p_{w}L_{w}+p_{w}^{*}L_{w}^{*}\right) ,$

\strut

for $p_{v},$ $p_{w},$ $p_{w}^{*}\in \Bbb{C}$. Here, $p_{w}^{*}$ is just a
complex number. Remark that $p_{w}^{*}$ is not a conjugate $\overline{p_{w}}$
of $p_{w}$ in $\Bbb{C}.$ The above expression of $a$ is said to be the
Fourier expansion of $a.$ Define the support $\Bbb{F}^{+}(G:a)$ of $a$ by

\strut

\begin{center}
$\Bbb{F}^{+}(G:a)=V(G:a)\cup FP(G:a)$
\end{center}

where

\begin{center}
$V(G:a)=\{v\in V(G):p_{v}\neq 0$ in (0.1)$\}$
\end{center}

and

\begin{center}
$
\begin{array}{ll}
FP(G:a)= & \{w\in FP(G):p_{w}\neq 0\text{ in (0.1)}\} \\ 
& \text{ \ }\cup \{w^{\prime }\in FP(G):p_{w}^{*}\neq 0\text{ in (0.1)}\}.
\end{array}
$
\end{center}

\strut

Notice that if $V(G:a)\neq \emptyset ,$ then $\underset{v\in V(G:a)}{\sum }%
p_{v}L_{v}$ is contained in the diagonal subalgebra $D_{G}.$ Thus we have
the canonical conditional expectation $E:W^{*}(G)\rightarrow D_{G},$ defined
by

\strut

\begin{center}
$E\left( a\right) =\underset{v\in V(G:a)}{\sum }p_{v}L_{v},$
\end{center}

\strut

for all $a$\ in $W^{*}(G)$ with its Fourier expansion (0.1). Then the
algebraic pair $\left( W^{*}(G),E\right) $ is a $W^{*}$-probability space
with amalgamation over $D_{G}$ (See [16]). This structure is called the
graph $W^{*}$-probability space over its diagonal subalgebra $D_{G},$ and
all elements in $(W^{*}(G),$ $E)$ are said to be $D_{G}$-valued random
variables. It is easy to check that the conditional expectation $E$ is
faithful in the sense that if $E(a^{*}a)=0_{D_{G}},$ for $a\in W^{*}(G),$
then $a=0_{D_{G}}.$

\strut

In [16] and [17], we computed the $D_{G}$-valued moments and cumulants of an
arbitrary $D_{G}$-valued random variable $a$ having its Fourier expansion
(0.1). In particular, by using the $D_{G}$-cumulant formula, we got the $%
D_{G}$-valued mixed cumulants of $D_{G}$-valued random variables $a_{1}$ and 
$a_{2}$ and then we could find the $D_{G}$-freeness characterization which
is so abstract to use. However, by using this characterization, we can
characterize the $D_{G}$-freeness of generators $L_{w}$'s ($w$ $\in $ $\Bbb{F%
}^{+}(G)$) $D_{G}$-freeness condition;

\strut

\ \ \ \ \ $L_{w_{1}}$ and $L_{w_{2}}$ are free over $D_{G}$ in $(W^{*}(G),E) 
$

\begin{center}
$\Longleftrightarrow w_{1}$ and $w_{2}$ in $F^{+}(G)$ are diagram-distinct,
\end{center}

\strut

in the sense that $w_{1}$ and $w_{2}$ have different diagram on the graph $%
G, $ graphically. Based on this $D_{G}$-freeness condition, in this paper,
we will observe the $D_{G}$-free structure of the given graph $W^{*}$%
-algebra $W^{*}(G).$

\strut \strut \strut \strut \strut

Define the subset $loop(G)$ of $FP(G),$ by the collection of all loop finite
paths in $FP(G).$ We will consider the family in $loop(G)$

\strut

\begin{center}
$\mathcal{F}=\{l_{j}\in loop(G):$ mutually diagram-distinct$\}_{j=1}^{N}$
\end{center}

\strut \strut

and construct the corresponding $D_{G}$-semicircular system

\strut

\begin{center}
$\mathcal{L}_{N}=\left\{ L_{l_{j}}+L_{l_{j}}^{*}:l_{j}\in F\right\} .$
\end{center}

\strut

Indeed, the elements in $\mathcal{L}_{N}$ are free from each other over $%
D_{G},$ by the diagram-distinctness of $\mathcal{F}$ and they are all $D_{G}$%
-semicircular, by [17]. So, the family $\mathcal{L}_{N}$ is the $D_{G}$%
-semicircular system in $\left( W^{*}(G),E\right) .$ We show that the free
product structure of the $W^{*}$-subalgebra $W^{*}(\mathcal{L}_{N},D_{G})$
generated by the $D_{G}$-semicircular system $\mathcal{L}_{N}$ and the
diagonal subalgebra $D_{G}$ ;

\strut

\begin{center}
$
\begin{array}{ll}
\left( W^{*}(\mathcal{L}_{N},D_{G}),E\right) & \,\simeq \left( W^{*}(%
\mathcal{L}_{N},D_{N})\otimes D_{G},E_{N}\otimes \mathbf{1}\right) \\ 
& 
\begin{array}{l}
\simeq \left( W^{*}(\mathcal{L}_{N},D_{N}),E_{N}\right) \otimes \left( D_{G},%
\mathbf{1}\right) \\ 
=\underset{j=1}{\overset{N}{*_{D_{G}}}}\left(
(W^{*}(\{L_{l_{j}}\},D_{N}),E_{N})\otimes \left( D_{G},\mathbf{1}\right)
\right)
\end{array}
\end{array}
$
\end{center}

\strut

where $D_{N}=\overline{\Bbb{C}[L_{v_{j}}:l_{j}=v_{j}l_{j}v_{j},\,l_{j}\in F]}%
^{w},$ $E_{N}$ $=$ $E_{D_{N}}^{D_{G}}$ $\circ $ $E,$ and $\mathbf{1}$ is the
identity map on $D_{G}.$

\strut

Also we define a new family

\strut

\begin{center}
$\mathcal{R}=\{L_{w_{j}},L_{w_{j}}^{*}:w_{j}\in loop^{c}(G)\}_{j=1}^{N},$
\end{center}

\strut

where $loop^{c}(G)=FP(G)\,\setminus \,loop(G).$ Then since the distinctness
of non-loop finite paths is the diagram-distinctness on them, $\{L_{w_{j}},$ 
$L_{w_{j}}^{*}\}$'s, for $j$ $=$ $1,$ $...,$ $N,$ in $\mathcal{R}$ are free
from each other over $D_{G}.$ Moreover, by [17], $L_{w_{j}}$ and $%
L_{w_{j}}^{*}$ are $D_{G}$-valued R-diagonal. We call this family $\mathcal{R%
},$ the $D_{G}$-valued R-diagonal system. Similar to the $D_{G}$%
-semicircular system case, we will define a $W^{*}$-subalgebra $W^{*}(%
\mathcal{R},D_{G})$ of the graph $W^{*}$-algebra $W^{*}(G),$ called the
R-diagonal subalgebra of $\left( W^{*}(G),E\right) .$ We can see that

\strut

\begin{center}
$\left( W^{*}(\mathcal{R},D_{G}),E_{\mathcal{R}}\right) =\left( \underset{j=1%
}{\overset{N}{\,\,\,*_{D_{\mathcal{R}}}}}\left( W^{*}(\{L_{w_{j}}\},D_{%
\mathcal{R}}),E_{D_{\mathcal{R}}}^{D_{G}}\circ E\right) \right) \otimes
\left( D_{G},\mathbf{1}\right) ,$
\end{center}

\strut \strut

where $D_{\mathcal{R}}=\overline{\Bbb{C}%
[\{L_{v_{1}},L_{v_{2}}:w_{j}=v_{1}w_{j}v_{2}\}]}^{w}$ and $E_{\mathcal{R}%
}=E\mid _{W^{*}(\mathcal{R},D_{G})}.$

\strut

We will also define the $D_{G}$-free building blocks of the graph $W^{*}$%
-algebra $W^{*}(G)$ and we prove that

\strut

\begin{center}
$
\begin{array}{ll}
\left( W^{*}(G),E\right) & \,=\left( D_{G},E\right) \\ 
& 
\begin{array}{l}
\\ 
\,\,\,\,*_{D_{G}}\left( \underset{l\in Loop(G)}{*_{D_{G}}}\left(
W^{*}(\{L_{l}\},D_{G}),E\right) \right)
\end{array}
\\ 
& 
\begin{array}{l}
\\ 
\,\,\,\,*_{D_{G}}\left( \underset{w\in loop^{c}(G)}{*_{D_{G}}}\left(
W^{*}(\{L_{w}\},D_{G}),E\right) \right) ,
\end{array}
\end{array}
$
\end{center}

\strut

where $Loop(G)$ is the set of all basic loops in $loop(G).$ A loop $l$ is
basic if there is no loop $w$ and a natural number $k$ $\in $ $\Bbb{N}$ $%
\setminus $ $\{1\}$ such that $l$ $=$ $w^{k}.$ This free product structure
of $\left( W^{*}(G),E\right) $ is nice for studying the subalgebras of $%
W^{*}(G).$

\strut

Finally, by considering the subalgebra inclusion, we can get the following
free product structure of $\left( W^{*}(G),E\right) ,$

\strut

\begin{center}
$\left( W^{*}(G),E\right) =\left( D_{G},E\right) *_{D_{G}}\left( \underset{%
e\in E(G)}{\,*_{D_{G}}}\left( W^{*}(\{L_{e}\},D_{G}),E\right) \right) .$
\end{center}

\strut \strut \strut \strut

\strut

\strut \strut \strut

\section{Semicircular System}

\bigskip

\subsection{The $D_{G}$-Semicircular System}

\bigskip

In this chapter, we will consider the amalgamated semicircular system
observed by Shlyaktenko (See [10]), in our graph structure. Throughout this
chapter, let $G$ be a countable directed graph and $(W^{*}(G),E)$, the graph 
$W^{*}$-probability space over the diagonal subalgebra $D_{G}.$

\bigskip \strut

\begin{definition}
Let $B$ be a von Neumann algebra and $A,$ a von Neumann algebra over $B$ and
let $F$ $:$ $A$ $\rightarrow $ $B$ be a conditional expectation. Then $(A,$ $%
F)$ is the amalgamated $W^{*}$-probability space over $B.$ The $B$-valued
random variable $x$ $\in $ $(A,$ $F)$ is said to be a $B$-semicircular
element if $x$ is self-adjoint and if the only nonvanishing $B$-cumulant of $%
x$ is the second $B$-cumulant of $x$. Let $x_{1},$ $...,$ $x_{s}$ be
self-adjoint $B$-valued random variables in $(A,$ $F),$ where $s$ $\in $ $%
\mathbb{N}.$ We say that the set $S$ $=$ $\{x_{1},$ $...,$ $x_{s}\}$ is a $B$%
\textbf{-semicircular family} if all $x_{j}$'s are $B$-semicircular, for $\ j
$ $=$ $1,$ $...,$ $s.$ The $B$-semicircular family $S$ is said to be a $B$%
\textbf{-semicircular system} if $x_{1},$ $...,$ $x_{s}$ are free from each
other over $B,$ in $(A,$ $F).$ The algebra generated by a $B$-semicircular
system and $B$ is called the $B$\textbf{-semicircular (sub)algebra} of $A.$
\end{definition}

\bigskip \strut \strut

Assume that we have a one-vertex directed graph $H.$ Then the diagonal
subalgebra $D_{H}$ $=$ $\Bbb{C}.$ So, in this case, the $D_{H}$%
-semicircularity is same as the Voiculescu's semicircularity. We will define
the lattice path model $LP_{n}^{*},$

\bigskip

\begin{center}
$LP_{n}^{*}=\{L:$ lattice path having the $*$-axis-property$\}$
\end{center}

\bigskip

(See [16] and [17]). Take $L\in LP_{n}^{*}.$ Then we have a (non-unique)
corresponding lattice path $l_{w_{1},...,w_{n}}^{u_{1},...,u_{n}}$ of the $%
D_{G}$-valued random variable $L_{w_{1}}^{u_{1}}$ ... $L_{w_{n}}^{u_{n}},$
where $u_{j}\in \{1,*\},$ in some graph $W^{*}$-probability space $%
(W^{*}(G), $ $E)$ over $D_{G}.$

\bigskip \strut \strut \strut

\begin{theorem}
(See [17]) The $D_{G}$-valued random variables $a_{l}=L_{l}+L_{l}^{*}$ are $%
D_{G}$-semicircular, for all $\ l$ $=$ $vlv$ $\in $ $loop(G),$ with $v$ $\in 
$ $V(G).$ In particular, we have that

\strut 

$\ \ \ \ \ \ \ \ \ \ \ \ \ k_{n}\left( a_{l},...,a_{l}\right) =\left\{ 
\begin{array}{ll}
2L_{v} & \text{if }n=2 \\ 
&  \\ 
0_{D_{G}} & \text{otherwise,}
\end{array}
\right. $

and

$\ \ \ \ \ \ \ \ \ \ \ \ \ E\left( a_{l}^{n}\right) =\left\{ 
\begin{array}{ll}
c_{\frac{n}{2}}\cdot \left( 2L_{v}\right) ^{\frac{n}{2}} & \text{if }n\text{
is even} \\ 
&  \\ 
0_{D_{G}} & \text{if }n\text{ is odd,}
\end{array}
\right. $

\strut 

for all $n\in \Bbb{N},$ where $c_{k}$ $=$ $\frac{1}{k+1}$ $\left( 
\begin{array}{l}
2k \\ 
\,\,k
\end{array}
\right) $ is the $k$-th Catalan number.$\ \square $
\end{theorem}

\bigskip \strut

We will consider the $D_{G}$-semicircular family $\{a_{1},$ $...,$ $a_{N}\},$
for $N$ $\in $ $\mathbb{N}.$ Let $l_{j}=v_{j}l_{j}v_{j}\in loop(G),$ for $%
j=1,...,N$, with $v_{j}\in V(G).$ Assume that loops $l_{1},$ $...,$ $l_{N}$
are mutually \textbf{diagram-distinct}. Define $D_{G}$-valued random
variables $a_{1},...,a_{N},$

\bigskip

\begin{center}
$a_{j}\overset{def}{=}L_{l_{j}}+L_{l_{j}}^{*},$ \ for all \ $j=1,...,N.$
\end{center}

\strut

Again, remark that we assumed that $l_{j}$'s are mutually diagram-distinct.
So, $a_{1},$ ..., $a_{N}$ are free from each other over $D_{G}$ in $%
(W^{*}(G),$ $E)$ and hence the $D_{G}$-semicircular family $\{a_{1},$ $...,$ 
$a_{N}\}$ is a $D_{G}$-semicircular system in $(W^{*}(G),$ $E).$ So, we have
the $D_{G}$-semicircular system, in $W^{*}(G),$ induced by the mutually
diagram-distinct loops in $FP(G).$ i.e., the set

\bigskip

\begin{center}
$\mathcal{L}_{N}=\{a_{j}:l_{j}$'s are$\,$diagram-distinct in $%
loop(G)\}_{j=1}^{N}$
\end{center}

\bigskip

is the $D_{G}$-semicircular system in $(W^{*}(G),E).$

\strut

\strut

\subsection{$D_{G}$-Semicircular Subalgebra of $\left( W^{*}(G),E\right) $}

\bigskip \strut

\strut \strut

Now, we will construct the $D_{G}$-semicircular algebra $W^{*}\left( 
\mathcal{L}_{N},D_{G}\right) ,$ as a $W^{*}$-subalgebra of the graph $W^{*}$%
-algebra, generated by $\mathcal{L}_{N}$ and $D_{G}.$ Let

\strut

\begin{center}
$\mathcal{F}=\{l_{j}\in loop(G):j=1,...,N\}$
\end{center}

\strut

be a collection of mutually diagram-distinct loops in $FP(G)$ and let

$\strut $

\begin{center}
$\mathcal{L}_{N}=\{a_{j}=L_{l_{j}}+L_{l_{j}}^{*}:l_{j}\in \mathcal{F}\}.$
\end{center}

\strut

Then the family $\mathcal{L}_{N}$ is a $D_{G}$-semicircular system in $%
\left( W^{*}(G),E\right) $ and the $W^{*}$-subalgebra $W^{*}\left( \mathcal{L%
}_{N},D_{G}\right) $ is the $D_{G}$-semicurcular subalgebra of $W^{*}(G).$
The $D_{G}$-semicircular subalgebra $W^{*}\left( \mathcal{L}%
_{N},D_{G}\right) $ have the following free product structure which is very
natural by the very definition.

\strut

\begin{lemma}
Let $(W^{*}(G),E)$ be a graph $W^{*}$-probability space over the diagonal
subalgebra $D_{G}$ and let

$\strut $

$\ \ \ \ \ \ \mathcal{L}_{N}=\{a_{j}=L_{l_{j}}+L_{l_{j}}^{*}:l_{j}$'s are
diagram-distinct in $loop(G)\}_{j=1}^{N}.$

\strut 

Then the $W^{*}$-subalgebra $W^{*}\left( \mathcal{L}_{N},D_{G}\right) $ of $%
W^{*}(G)$ is a $D_{G}$-semicircular algebra satisfies that

\strut 

$\ \ \ \ \ \ \ \ \ \ \ \ \ \ \ W^{*}\left( \mathcal{L}_{N},D_{G}\right)
\simeq \underset{j=1}{\overset{N}{\text{ }*_{D_{G}}}}~W^{*}(a_{j},D_{G}).$

$\square $
\end{lemma}

\bigskip \strut

Let $\mathcal{L}_{N}$ be give as above. Assume that $l_{j}=v_{j}l_{j}v_{j},$
for $v_{j}$ $\in $ $V(G).$ (It is possible that $v_{i}$ $=$ $v_{k},$ for
some $i,$ $k$ in $\{1,$ $...,$ $N\}.$) Define the subalgebra $D_{N}$ of the
diagonal subalgebra $D_{G}$ by

\strut

\begin{center}
$D_{N}=\overline{\Bbb{C}[\{L_{v_{j}}:j=1,...,N\}]}^{w}.$
\end{center}

\strut

Trivially, $D_{N}\leq D_{G},$ as von Neumann algebras.

\strut

\begin{proposition}
(Also See [18]) Let $\mathcal{L}_{N}$ be the given $D_{G}$-semicircular
system and let $D_{N}$ be defined as above. As amalgamated $W^{*}$%
-probability spaces,

\strut 

$\ \ \ \ \ \ 
\begin{array}{ll}
\left( W^{*}\left( \mathcal{L}_{N},D_{G}\right) ,E\right)  & \simeq \left(
W^{*}\left( \mathcal{L}_{N},D_{N}\right) \otimes D_{G},E_{N}\otimes \mathbf{1%
}\right)  \\ 
& \simeq \left( W^{*}(\mathcal{L}_{N}),E_{N}\right) \otimes \left( D_{G},%
\mathbf{1}\right) ,
\end{array}
$

\strut 

where $E_{N}:W^{*}\left( \mathcal{L}_{N}\right) \rightarrow D_{N}$ is the
conditional expectation defined by $E_{N}$ $=$ $E_{D_{N}}^{D_{G}}$ $\circ $ $%
E$ and $\mathbf{1}(d)$ $=$ $d,$ $\forall $ $d$ $\in $ $D_{G}.$
\end{proposition}

\strut \strut \strut \strut \strut \strut

\begin{proof}
As $W^{*}$-algebras, $W^{*}(\mathcal{L}_{N},D_{G})\simeq W^{*}(\mathcal{L}%
_{N},D_{N})\otimes D_{G}.$ Indeed, without loss of generality, take $a\in
W^{*}(\mathcal{L}_{N},D_{G})$ by

$\strut $

$\ \ \ \ \ \ \ \ \
a=d_{1}a_{l_{i_{1}}}^{k_{1}}d_{2}a_{l_{i_{2}}}^{k_{2}}...d_{n}a_{l_{i_{n}}}^{k_{n}} 
$ \ and \ $a_{l_{j}}=L_{l_{j}}+L_{l_{j}}^{*}$

\strut

where $d_{1},...,d_{n}\in D_{G},$ $k_{1},...,k_{n}\in \Bbb{N}$ and $%
(i_{1},...,i_{n})$ $\in $ $\{1,$ $...,$ $N\}^{n},$ $n\in \mathbb{N}.$
Observe that, for any $j$ $\in $ $\{1,$ $...,$ $N\},$ we have that

$\strut $

$\ \ \ \ \ \ \ \ \ a_{l_{j}}^{k}=\left( L_{l_{j}}+L_{l_{j}}^{*}\right)
^{k}=L_{l_{j}^{k}}+L_{l_{j}^{k}}^{*}+Q\left( L_{l_{j}},L_{l_{j}}^{*}\right)
, $

\strut

where $Q\in \mathbb{C}[z_{1},z_{2}]$. Also, observe that $%
L_{l_{j}}^{k_{1}}L_{l_{j}}^{*\,\,k_{2}},$ for any $k_{1}$ $,k_{2}$ $\in $ $%
\mathbb{N},$ satisfies that

\strut

$\ \ \
L_{l_{j}}^{k_{1}}L_{l_{j}}^{*\,%
\,k_{2}}=L_{l_{j}^{k_{1}}}L_{l_{j}^{k_{2}}}^{*}=\left\{ 
\begin{array}{ll}
L_{l_{j}^{k_{1}-k_{2}}}=L_{v_{j}}L_{l_{j}^{k_{1}-k_{2}}} & \text{if }%
k_{1}>k_{2} \\ 
L_{l_{j}^{k_{2}-k_{1}}}^{*}=L_{v_{j}}L_{l_{j}^{k_{2}-k_{1}}}^{*} & \text{if }%
k_{1}<k_{2} \\ 
L_{v_{j}} & \text{if }k_{1}=k_{2},
\end{array}
\right. $

and similarly,

\strut

$\ \ \
L_{l_{j}}^{*\,\,k_{1}}L_{l_{j}}^{%
\,k_{2}}=L_{l_{j}^{k_{1}}}^{*}L_{l_{j}^{k_{2}}}=\left\{ 
\begin{array}{ll}
L_{l_{j}^{k_{1}-k_{2}}}^{*}=L_{v_{j}}L_{l_{j}^{k_{1}-k_{2}}}^{*} & \text{if }%
k_{1}>k_{2} \\ 
L_{l_{j}^{k_{2}-k_{1}}}=L_{v_{j}}L_{l_{j}^{k_{2}-k_{1}}} & \text{if }%
k_{1}<k_{2} \\ 
L_{v_{j}} & \text{if }k_{1}=k_{2}.
\end{array}
\right. $

So,

\strut

$\ \ \ \ \ \ \ \ \ \ \ \ \ \ \ Q(L_{l_{j}},L_{l_{j}}^{*})=L_{v_{j}}\left(
Q(L_{l_{j}},L_{l_{j}}^{*})\right) L_{v_{j}},$

\strut

for all \ $j=1,...,N.$ Thus

\strut

(1.1) $\ \ \ \ \ \ \
a_{l_{j}}^{k}=L_{v_{j}}a_{l_{j}}^{k}=L_{v_{j}}a_{l_{j}}^{k}L_{v_{j}},$ for
all $j=1,...,N.$

\strut \strut \strut

Now, consider that

\strut

$\ \ \ \ \ \ \ \ \ \ \ \ \ d_{j}=d_{j}^{N}+d_{j}^{\prime },$ $\forall
j=1,...,N.$

\strut

where $d_{j}^{N}=\sum_{j=1}^{N}L_{v_{j}}d_{j}L_{v_{j}}$ and $d_{j}^{\prime
}=d_{j}-d_{j}^{N}$ in $D_{G}.$ So, we can rewrite that

\strut

$\ \ a=\left( d_{1}^{N}+d_{1}^{\prime }\right) a_{l_{i_{1}}}^{k_{1}}\left(
d_{2}^{N}+d_{2}^{\prime }\right) a_{l_{i_{2}}}^{k_{2}}...\left(
d_{n}^{N}+d_{n}^{\prime }\right) a_{l_{i_{n}}}^{k_{n}}$

\strut

$\ \ \ \
=d_{1}^{N}a_{l_{i_{1}}}^{k_{1}}d_{2}^{N}a_{l_{i_{2}}}^{k_{2}}...d_{n}^{N}a_{l_{i_{n}}}^{k_{n}}+d_{1}^{\prime }a_{l_{i_{1}}}^{k_{1}}d_{2}^{\prime }a_{l_{i_{2}}}^{k_{2}}...d_{n}^{\prime }a_{l_{i_{n}}}^{k_{n}} 
$

\strut

$\ \ \ \
=d_{1}^{N}a_{l_{i_{1}}}^{k_{1}}d_{2}^{N}a_{l_{i_{2}}}^{k_{2}}...d_{n}^{N}a_{l_{i_{n}}}^{k_{n}}, 
$

\strut \strut

by (1.1). This shows that $a=a\otimes 1\in W^{*}(\mathcal{L}%
_{N},D_{N})\otimes 1$ and

\strut

$\ \ \ \ \ \ \ \ E\left( a\right) =E_{D_{N}}^{D_{G}}\circ
E(a)=E_{N}(a)=E_{N}\otimes \mathbf{1}(a\otimes 1).$

\strut

Trivially, if $a\in D_{G}\subset W^{*}(\mathcal{L}_{N},D_{G}),$ then $%
a=1\otimes a\in 1\otimes D_{G}.$ Futhermore, if $a\in D_{G}$ in $W^{*}(%
\mathcal{L}_{N},D_{G}),$ then

\strut

\ $\ \ \ \ \ \ \ \ \ \ E(a)=a=1\otimes a=E_{N}\otimes \mathbf{1}(1\otimes
a). $

\strut

Now, consider $W^{*}(\mathcal{L}_{N},D_{N}).$ By the previous lemma,
similarly, we have that

\strut

(1.2) $\ \ \ \ W^{*}(\mathcal{L}_{N},D_{N})=W^{*}(\{a_{1}%
\},D_{N})*_{D_{N}}...*_{D_{N}}W^{*}(\{a_{N}\},D_{N}).$

\strut

Indeed, the $D_{G}$-semicircular elements $a_{i}$ and $a_{j}$ in $L_{N}$ are
free over $D_{N}$ in $W^{*}(L_{N},D_{N}).$ Clearly, since $D_{N}\subset
D_{G} $ and since $a_{i}$ and $a_{j}$ are free over $D_{G},$ they are free
over $D_{N}.$ Therefore, the formula (3.1.2) holds true with respect to the
(compressed) conditional expectation

\strut

\ \ \ \ \ \ \ \ \ \ \ \ \ \ \ \ \ \ \ \ \ \ \ $E_{N}=E\mid
_{W^{*}(L_{N},D_{N})}=E_{D_{N}}^{D_{G}}\circ E$,

\strut

on $W^{*}(\mathcal{L}_{N},D_{N})$.
\end{proof}

\strut \strut

\begin{corollary}
Let $\mathcal{F}=\{l_{j}:l_{j}=v_{0}l_{j}v_{0}\}_{j=1}^{N}$ be a family of
mutually diagram-distinct loops in $loop(G),$ where $v_{0}$ $\in $ $V(G)$.
If the collection

\strut 

$\ \ \ \ \ \ \ \ \ \ \ \ \ \ \ \ \ \ \mathcal{L}_{N}=\{\frac{1}{\sqrt{2}}%
a_{l_{j}}:l_{j}\in \mathcal{F}\}$

\strut 

is a $D_{G}$-semicircular system induced by the family $\mathcal{F},$ then

\strut 

$\ \ \ \ \ \ \ \ \ \ \ \left( W^{*}\left( \mathcal{L}_{N},D_{G}\right)
,E\right) \simeq \left( L(F_{N}),tr\right) \otimes \left( D_{G},\,\mathbf{1}%
\right) ,$

\strut 

in the sense of Voiculescu, where $tr$ is the canonical $II_{1}$-trace of
the free group factor $L(F_{K}),$ $\forall $ $K$ $\in $ $\Bbb{N}$, and where 
$\mathbf{1}(d)$ $=$ $d,$ $\forall $ $d$ $\in $ $D_{G}.$ $\square $
\end{corollary}

\strut \strut \strut

By the previous proposition and corollary, we can have the following fact ;

\strut

\begin{theorem}
Let $\mathcal{F}=\{l_{k1},...,l_{kn_{k}}:l_{kj}=v_{k}l_{kj}v_{k},$ $%
j=1,...,n_{k}\}_{k=1}^{N}$ be the collection of $\sum_{k=1}^{N}n_{k}$%
-mutually diagram-distinct loops in $FP(G).$ Assume that $v_{k_{1}}\neq
v_{k_{2}},$ for all pair $(k_{1},k_{2})\in \Bbb{N}^{2}.$ If

\strut 

\ $\ \ \ \ \ \ \ \mathcal{L}=\{\frac{1}{\sqrt{2}}a_{l_{k_{j}}}:j=1,...,n_{k}%
\}_{k=1}^{N}$

and

\ $\ \ \ \ \ \ \ D_{\mathcal{L}}=\overline{\Bbb{C}[\{L_{v_{k}}:k=1,...,N\}]}%
^{w}$

\strut 

are the corresponding $D_{G}$-semicircular system and the subalgebra of $%
D_{G}$ by $L$, respectively, and

\strut 

\ \ $\ \ \ \ \ \ \ \ \ \ \ \ \ \ \ \ E_{\mathcal{L}}=E_{D_{L}}^{D_{G}}\circ
E,$

then

\strut 

$\ \ \ \left( W^{*}(\mathcal{L},D_{G}),E\right) \simeq \left( \underset{k=1}{%
\overset{N}{\,\,\,\,*_{D_{\mathcal{L}}}}}\left( \left(
L(F_{n_{k}}),tr\right) \otimes \left( D_{\mathcal{L}},\mathbf{1}\right)
\right) \right) \otimes \left( D_{G},\mathbf{1}\right) .$

$\square $
\end{theorem}

\strut \strut \strut \strut \strut

\begin{corollary}
Let $\mathcal{F}_{1}=\{l_{1}^{1},...,l_{N}^{1}:l_{j}^{1}=v_{1}l_{j}^{1}v_{1}%
\}$ and $\mathcal{F}_{2}=%
\{l_{1}^{2},...,l_{N}^{2}:l_{j}^{2}=v_{2}l_{j}^{2}v_{2}\}$ be the $N$%
-mutually diagram-distinct families of loops in $FP(G),$ where $v_{1}\neq
v_{2}\in V(G)$ are fixed. Let

\strut \strut 

$\ \ \ \ \ \ \ \mathcal{L}_{1}=\left\{ a_{k}=\frac{1}{\sqrt{2}}\left(
L_{l_{k}^{1}}+L_{l_{k}^{1}}^{*}\right) :k=1,...,N\right\} $

and

$\ \ \ \ \ \ \ \mathcal{L}_{2}=\left\{ b_{k}=\frac{1}{\sqrt{2}}\left(
L_{l_{k}^{2}}+L_{l_{k}^{2}}^{*}\right) :k=1,...,N\right\} $

\strut 

be the corresponding $D_{G}$-semicircular systems, respectively. Then two $%
D_{G}$-semicircular subalgebras $\left( W^{*}(\mathcal{L}_{1},D_{G}),E%
\right) $ and $\left( W^{*}(\mathcal{L}_{2},D_{G}),E\right) $ are free over $%
D_{G}$ and they are isomorphic, as amalgamated $W^{*}$-probability spaces. $%
\square $
\end{corollary}

\strut \strut \strut \strut \strut

The above corollary shows us how to construct the isomorphic semicircular
subalgebras in the graph $W^{*}$-probability space from two vertices having
the same number of loops. (Assume that the vertex $v_{1}$ has its loops $%
l_{1}^{1},$ $...,$ $l_{n_{1}}^{1}$ and $v_{2}$ has its loops $l_{1}^{2},$ $%
...,$ $l_{n_{2}}^{2}$ and suppose that $n_{1}<n_{2}.$ Then we can choose $%
n_{1}$-loops of $v_{2},$ \ $l_{i_{1}}^{2},$ $...,$ $l_{i_{n_{1}}}^{2}.$ And
then we can apply the above corollary for them.)\strut

\strut

One of the most interesting example of $D_{G}$-semicircular subalgebra is as
follows ;

\strut \strut \strut \strut \strut

\begin{example}
Let $G$ be a directed graph with

\strut 

$\ \ \ \ \ \ \ V(G)=\{v\}$ and $E(G)=\{l_{1},...,l_{N}:l_{j}=vl_{j}v\}.$

\strut 

Note that $D_{G}=\Delta _{1}=\Bbb{C}.$ Also, note that the projection $L_{v}$
$=$ $1_{D_{G}}$ $=$ $1_{\Bbb{C}}$ $=$ $1$ $\in $ $\Bbb{C}.$ We have the
graph $W^{*}$-probability space $\left( W^{*}(G),E\right) $ over $\Bbb{C}$.
Consider the $W^{*}$-subalgebra $W^{*}\left( \mathcal{L}_{N},D_{G}\right) $ $%
=$ $W^{*}(\mathcal{L}_{N}),$ where

\strut 

$\ \ \ \ \ \ \ \ \ \ \ \mathcal{L}_{N}=\left\{ \frac{1}{\sqrt{2}}%
a_{j}:\,a_{j}=L_{l_{j}^{n_{j}}}+L_{l_{j}^{n_{j}}}^{*}\right\} _{j=1}^{N},$

\strut 

where $n_{j}\in \Bbb{N},$ $j=1,...,N.$ Then $\mathcal{L}_{N}$ is a $D_{G}$%
-semicircular system, too. Definitely it is a $D_{G}$-semicircular family.
Since $\mathcal{F}$ $=$ $\{l_{j}^{n_{j}}$ $:$ $j$ $=$ $1,$ $...,$ $N\}$ is
consists of mutually diagram-distinct loops in $FP(G),$ $a_{j}$'s are free
from each other and hence $\mathcal{L}_{N}$ is a $D_{G}$-semicircular
system. We have that

\strut 

$\ k_{2}\left( \frac{1}{\sqrt{2}}a_{j},\,\frac{1}{\sqrt{2}}a_{j}\right) =%
\frac{1}{2}k_{2}\left( a_{j},a_{j}\right) $

\strut 

$\ \ \ \ \ =\frac{1}{2}\,k_{2}\left(
L_{l_{j}^{n_{j}}}+L_{l_{j}^{n_{j}}}^{*},%
\,L_{l_{j}^{n_{j}}}+L_{l_{j}^{n_{j}}}^{*}\right) $

\strut 

$\ \ \ \ \ =\frac{1}{2}\,\underset{(r_{1},r_{2})\in \{1,*\}^{2}}{\sum }%
k_{2}\left( L_{l_{j}^{n_{j}}}^{r_{1}},\,L_{l_{j}^{n_{j}}}^{r_{2}}\right) $

\strut 

$\ \ \ \ \ =\frac{1}{2}\left( \mu _{l_{j}^{n_{j}},l_{j}^{n_{j}}}^{1,*}\cdot
E\left( L_{l_{j}^{n_{j}}},\,L_{l_{j}^{n_{j}}}^{*}\right) +\mu
_{l_{j}^{n_{j}},l_{j}^{n_{j}}}^{*,1}\cdot E\left(
L_{l_{j}^{n_{j}}}^{*},\,L_{l_{j}^{n_{j}}}\right) \right) $

\strut 

$\ \ \ \ \ =\frac{1}{2}\left( L_{v}+L_{v}\right) =\frac{1}{2}\cdot 2=1,$

\strut 

for all \ $j=1,...,N.$ So, we can get that

\strut 

$\ \ \ \ \ \ 
\begin{array}{ll}
\left( W^{*}(L_{N}),E\right)  & =\left( W^{*}\left( \mathcal{L}%
_{N},D_{G}\right) ,E\right)  \\ 
& =\left( W^{*}\left( \mathcal{L}_{N},D_{N}\right) ,E_{N}\right) \otimes
\left( D_{G},\mathbf{1}\right)  \\ 
& =\left( W^{*}(\mathcal{L}_{N}),E_{N}\right) \otimes \left( \Bbb{C},\mathbf{%
1}\right)  \\ 
& =\left( L(F_{N}),tr\right) ,
\end{array}
$

\strut 

where $tr:L(F_{N})\rightarrow \Bbb{C}$ is the canonical $II_{1}$-trace on
the free group factor $L(F_{K}),$ $\forall K\in \Bbb{N}.$ So, the graph $%
W^{*}$-probability space $\left( W^{*}(G),E\right) $ contains the free group
factor $L(F_{N})$ which is isomorphic to the $D_{G}$-semicircular subalgebra 
$W^{*}(\mathcal{L}_{N}),$ generated by $\mathcal{L}_{N}.$
\end{example}

\strut \strut

\strut

\strut \strut

\section{R-diagonal Systems}

\strut

\strut

In this chapter, similar to Chapter 2, we will consider the special $W^{*}$%
-subalgebra of the graph $W^{*}$-probability space $\left( W^{*}(G),E\right) 
$ over the diagonal subalgebra $D_{G}.$ As we defined the $D_{G}$%
-semicircular systems in $W^{*}(G),$ we will define the ($D_{G}$-valued)
R-diagonal systems in $W^{*}(G).$ Recall that if $w\in \,loop^{c}(G)$ is a
(non-loop) finite path in $\Bbb{F}^{+}(G),$ then the $D_{G}$-valued random
variables $L_{w}$ and $L_{w}^{*}$ are R-diagonal over $D_{G}$ (See [17]).
Take a finite family

\strut

\begin{center}
$\mathcal{F}=\{w_{j}:w_{j}\in loop^{c}(G)\}_{j=1}^{N},$ $N\in \Bbb{N}$.
\end{center}

\strut

Define a ($D_{G}$-valued) R-diagonal family induced by $\mathcal{F}$ by

\strut

\begin{center}
$R=\{L_{w},L_{w}^{*}:w\in \mathcal{F}\}.$
\end{center}

\strut

Notice that, by the (diagram-)distinctness of $w_{j}$'s in $\mathcal{F}$,
the subfamilies $\{L_{w_{1}},L_{w_{1}}^{*}\},$ ..., $%
\{L_{w_{N}},L_{w_{N}}^{*}\}$ are free from each other over $D_{G}$ in $%
\left( W^{*}(G),E\right) .$ We will observe that the R-diagonal subalgebra $%
W^{*}(R,D_{G})$ satisfies that

\strut

\begin{center}
$\left( W^{*}(R,D_{G}),E\right) =\underset{j=1}{\overset{N}{\,\,\,*_{D_{R}}}}%
\left( \left( W^{*}(\{L_{w_{j}}\},D_{R}),E_{R}\right) \otimes \left( D_{G},%
\mathbf{1}\right) \right) ,$
\end{center}

\strut

where $D_{R}$ is the subalgebra generated by the projections

\strut

\begin{center}
$\{L_{v_{1}},L_{v_{2}}:w=v_{1}wv_{2},\forall w\in \mathcal{F}\}$.
\end{center}

\strut

and $E_{R}=E\mid _{W^{*}(\{L_{w_{j}}\}_{j=1}^{N},\,D_{R})}$ is the
restricted conditional expectation onto $D_{R}.$

\strut

\strut

\strut

\subsection{$D_{G}$-valued R-diagonal Systems}

\strut

\strut

\strut

Let $G$ be a countable directed graph and $\left( W^{*}(G),E\right) $, the
graph $W^{*}$-probability space over the diagonal subalgebra $D_{G}.$
Throughout this section, we will fix the following finite family,

\strut

\begin{center}
$\mathcal{F}=\{w_{j}:w_{j}\,\in loop^{c}(G)\}_{j=1}^{N}\subset FP(G),$
\end{center}

\strut

where $N\in \Bbb{N}.$ Notice that all elements in the family $\mathcal{F}$
are non-loop finite paths and hence the corresponding $D_{G}$-valued random
variables are $D_{G}$-valued R-diagonal elements (See [17]).

\strut

\begin{definition}
Let $a\in \left( W^{*}(G),E\right) $ be a $D_{G}$-valued random variable.
The $D_{G}$-valued random variable is said to be an ($D_{G}$-valued)
R-diagonal element if it has the only nonvanishing $D_{G}$-valued cumulants
having their forms of

\strut 

$\ \ \ \ \ \ \ \ \ k_{n}\left( a,a^{*},...,a,a^{*}\right) $ \ and \ $%
k_{n}\left( a^{*},a,...,a^{*},a\right) ,$

\strut 

for all $n\in 2\Bbb{N}$ \ (See [19] and [17]). Clearly, if $a$ is a
R-diagonal, then automatically $a^{*}$ is R-diagonal. (In other words, the
pair $\left( a,\text{ }a^{*}\right) $ is a $D_{G}$-valued R-diagonal pair.
Also, see [19].) Suppose we have a collection

\strut 

$\ \ \ \ \ \ \ \ R=\{a_{1},...,a_{N}:a_{j}$ is R-diagonal over $D_{G}\}.$

\strut 

Then the family $R$ is called the ($D_{G}$-valued) R-diagonal system. The
subalgebra $W^{*}(R,D_{G})$ is called the ($D_{G}$-valued) R-diagonal
subalgebra, induced by $R$ in $\left( W^{*}(G),E\right) .$
\end{definition}

\strut \strut \strut \strut

In [17], we showed the following theorem;

\strut

\begin{theorem}
Let $w\in FP(G).$ Then the $D_{G}$-valued random variables $L_{w}$ and $%
L_{w}^{*}$ are R-diagonal over $D_{G}$ in $\left( W^{*}(G),E\right) .$ $%
\square $
\end{theorem}

\strut \strut \strut

Notice that all $D_{G}$-semicircular elements are R-diagonal, by definition.
But we will restrict our interests to the R-diagonal systems consisting of
the non-loop finite paths. Recall that if $w_{1}$ $\neq $ $w_{2}$ $\in $ $%
\,loop^{c}(G)$, then the $D_{G}$-valued R-diagonal elements $L_{w_{1}}$ and $%
L_{w_{2}}$ are free over $D_{G},$ in $\left( W^{*}(G),E\right) ,$ since the
distinctness of non-loop finite paths is equivalent to the
diagram-distinctness of them.

\strut

\strut

\strut

\strut

\subsection{$D_{G}$-valued R-diagonal Subalgebras}

\strut

\strut

In this section, we will consider the $D_{G}$-valued R-diagonal subalgebras
of $W^{*}(G),$ generated by the fixed R-diagonal systems and $D_{G}.$ As
before, let $\mathcal{F}$ be a finite family consisting of $N$-mutually
diagram-distinct non-loop finite paths and let $R=\{L_{w}:w\in \mathcal{F}%
\}. $ As we saw in the previous section, the family $R$ is the R-diagonal
system in $\left( W^{*}(G),E\right) .$ Therefore, we can get the following
result ;

\strut \strut

\begin{proposition}
Let $\mathcal{F}$ and $R$ be given as before and let $W^{*}(R,D_{G})$ be the
R-diagonal subalgebra of $W^{*}(G).$ Then

\strut \strut 

$\ \ \ \ \ \ \ \left( W^{*}(R,D_{G}),E_{R}\right) =\underset{j=1}{\overset{N%
}{\,\,\,*_{D_{G}}}}\left( W^{*}(\{L_{w_{j}}\},D_{G}),E_{j}\right) ,$

\strut 

where $E_{R}=E\mid _{W^{*}(R,D_{G})}$ and $E_{j}=E\mid
_{W^{*}(\{L_{w_{j}}\},D_{G})},$ for all \ $j=1,...,N.$ $\square $
\end{proposition}

\strut \strut \strut \strut

Similar to the previous chapter, we will define

\strut

\begin{center}
$D_{R}\overset{def}{=}\overline{\Bbb{C}[\{L_{v_{1}},L_{v_{2}}:w=v_{1}wv_{2},%
\,w\in \mathcal{F\}}]}^{w}.$
\end{center}

\strut

Notice that, for the given R-diagonal system $R,$ the von Neumann algebra $%
D_{R}$ should not be $\Bbb{C}$, because $\mathcal{F}$ is consists of all
non-loop finite paths which are mutually distinct. For example, if the
family $\mathcal{F}$ $=$ $\{w_{0}$ $=$ $v_{1}$ $w_{0}$ $v_{2}\},$ where $%
v_{1}$ $\neq $ $v_{2}$ in $V(G),$ and $R$ $=$ $\{L_{w_{0}}$ $:$ $w_{0}$ $\in 
$ $\mathcal{F}\}.$ Then

\strut

\begin{center}
$D_{R}=\overline{\Bbb{C}[\{L_{v_{1}},L_{v_{2}}\}]}^{w}\simeq \Delta _{2},$
\end{center}

\strut

where $\Delta _{2}$ is the subalgebra of the matricial algebra $M_{2}(\Bbb{C}%
)$ generated by all diagonal matrices. Notice that, for the inclusion $%
D_{R}\subset D_{G},$ there exists the well-determined canonical conditional
expectation $E_{D_{R}}^{D_{G}}$ $:$ $D_{G}$ $\rightarrow $ $D_{R}.$ Also,
notice that $D_{R}$ $\subset $ $W^{*}(R).$

\strut

\begin{proposition}
Let $\mathcal{F}$ and $R$ be given as before and let $W^{*}(R,D_{G})$ is the
R-diagonal subalgebra of $W^{*}(G).$ Then

\strut 

$\ \ \ \left( W^{*}(R,D_{G}),E_{R}\right) =\left(
W^{*}(R,D_{R}),E_{D_{R}}^{D_{G}}\circ E\right) \otimes \left( D_{G},\mathbf{1%
}\right) ,$

\strut 

where $E_{R}=E\mid _{W^{*}(R,D_{G})}$ and $\mathbf{1}$ is the identity map
on $D_{G}.$
\end{proposition}

\strut

\begin{proof}
Let $a\in \left( W^{*}(R,D_{G}),E_{R}\right) $ be the nonzero $D_{G}$-valued
random variable such that

\strut

$\ \ \ \ \ \ \ \ \
a=d_{1}L_{w_{i_{1}}}^{r_{1}}d_{2}L_{w_{2}}^{r_{2}}...d_{n}L_{w_{n}}^{r_{n}}$
\ \ or \ $a\in D_{G},$

\strut

where $d_{1},...,d_{n}\in D_{G},$ $r_{1},...,r_{n}\in \{1,*\}$ and $%
(i_{1},...,i_{n})\in \{1,...,N\}^{n},$ $n\in \Bbb{N}.$

\strut

Let $d_{j}=\underset{v_{j}\in V(G:d_{j})}{\sum }q_{v_{j}}L_{v_{j}},$ for all
\ $j=1,...,n.$ Then

\strut

$\ a=\underset{(v_{1},...,v_{n})\in \Pi _{j=1}^{n}V(G:d_{j})}{\sum }\,\left(
\Pi _{j=1}^{n}q_{v_{j}}\right) \left(
L_{v_{1}}L_{w_{i_{1}}}^{r_{1}}L_{v_{2}}L_{w_{2}}^{r_{2}}...L_{v_{n}}L_{w_{n}}^{r_{n}}\right) 
$

\strut

$\ \ \ =\underset{(v_{1},...,v_{n})\in \Pi _{j=1}^{n}V(G:d_{j})}{\sum }%
\,\left( \Pi _{j=1}^{n}q_{v_{j}}\right) \left( \Pi _{j=1}^{n}\delta
_{(v_{j},x_{j},y_{j}:r_{j})}\right) $

\strut

$\ \ \ \ \ \ \ \ \ \ \ \ \ \ \ \ \ \ \ \ \ \ \ \ \ \ \ \ \ \ \ \ \ \ \ \ \ \
\ \ \ \ \ \left(
L_{w_{i_{1}}}^{r_{1}}L_{w_{2}}^{r_{2}}...L_{w_{n}}^{r_{n}}\right) ,$

\strut

where $\delta _{(v_{j},x_{j},y_{j}:r_{j})}$ $=$ $\delta _{v_{j},x_{j}}$ if $%
r_{j}$ $=$ $1$ and $\delta _{(v_{j},x_{j},y_{j}:r_{j})}$ $=$ $\delta
_{v_{j},y_{j}}$ if $r_{j}$ $=$ $*$. Now, assume that there exists at least
one $j$ $\in $ $\{1,$ $...,$ $n\}$ such that $v_{j}$ $\notin $ $D_{R}.$ Then 
$\delta _{(v_{j},x_{j},y_{j}:r_{j})}$ $=$ $0,$ where either $r_{j}$ $=$ $1$
or $r_{j}$ $=$ $*.$ So, to make $a$ be nonzero, the elements $d_{1},$ $...,$ 
$d_{n}$ should be chosen in $D_{R}$ $\subset \,D_{G}.$ This shows that the
arbitrary element $x$ in the R-diagonal subalgebra $W^{*}\left(
R,D_{G}\right) ,$ we have that

\strut

$\ \ \ \ \ \ \ \ \ \ \ \ \ \ x=x\otimes 1\in W^{*}(R,D_{R})\otimes 1$

\strut

\ \ \ \ \ \ \ or $\ \ \ x=1\otimes a\in 1\otimes D_{G}$

\strut

\ \ \ \ \ \ \ or $\ \ \ x=x_{1}\otimes x_{2}\in W^{*}(R,D_{R})\otimes D_{G}.$

\strut

Also, if $%
a=d_{1}L_{w_{i_{1}}}^{r_{1}}d_{2}L_{w_{2}}^{r_{2}}...d_{n}L_{w_{n}}^{r_{n}}%
\in W^{*}(R,D_{R}),$ in $W^{*}(R,D_{G}),$ for $d_{1},$ $...,$ $d_{n}$ $\in $ 
$D_{R},$ then

\strut

$\ \ 
\begin{array}{ll}
E_{R}\left(
d_{1}L_{w_{i_{1}}}^{r_{1}}d_{2}L_{w_{2}}^{r_{2}}...d_{n}L_{w_{n}}^{r_{n}}%
\right) & =E_{D_{R}}^{D_{G}}\circ E\left(
d_{1}L_{w_{i_{1}}}^{r_{1}}d_{2}L_{w_{2}}^{r_{2}}...d_{n}L_{w_{n}}^{r_{n}}%
\right) \\ 
&  \\ 
& =E_{D_{R}}^{D_{G}}\circ E\left(
d_{1}L_{w_{i_{1}}}^{r_{1}}d_{2}L_{w_{2}}^{r_{2}}...d_{n}L_{w_{n}}^{r_{n}}%
\right) .
\end{array}
$

\strut

Trivially, if $a\in D_{G}$ in $W^{*}(R,D_{G}),$ then

\strut

$\ \ \ \ \ \ \ \ \ \ \ \ \ \ \ E_{R}(a)=E(a)=a=\mathbf{1}(a)$,

\strut

where $\mathbf{1}$ is the identity map on $D_{G}.$
\end{proof}

\strut

By the above proposition, we can conclude this section by the following
theorem ;

\strut

\begin{theorem}
Let $\mathcal{F}$ and $R$ be given as before and let $W^{*}(R,D_{G})$ be the
R-diagonal subalgebra of $W^{*}(G).$ Then

\strut 

$\ \ \left( W^{*}(R,D_{G}),E_{R}\right) =\left( \underset{j=1}{\overset{N}{%
\,\,\,*_{D_{R}}}}\left( W^{*}(\{L_{w_{j}}\},D_{R}),E_{D_{R}}^{D_{G}}\circ
E\right) \right) \otimes \left( D_{G},\mathbf{1}\right) .$
\end{theorem}

\strut

\begin{proof}
Notice that if two $D_{G}$-valued random variables $x$ and $y$ are free over 
$D_{G}$ in $\left( W^{*}(G),E\right) ,$ then they are free over $D_{R},$
since $D_{R}\subset D_{G}.$ i.e.e, since all mixed $D_{G}$-valued cumulants
of $p(x,x^{*})$ and $q(y,y^{*})$ vanish, for all $p,q\in \Bbb{C}%
[z_{1},z_{2}],$ all their mixed $D_{R}$-valued cumulants vanish, again.
Thus, we have that

\strut

$\ \ \ \ \ \ \left( W^{*}(R,D_{R}),E_{D_{R}}^{D_{G}}\circ E\right) =%
\underset{j=1}{\overset{N}{\,\,\,*_{D_{G}}}}\left(
W^{*}(\{L_{w_{j}}\},D_{G}),E_{D_{R}}^{D_{G}}\circ E\right) .$

\strut \strut

Therefore, by the previous proposition, we can get that

\strut \strut

$\ \ \ \left( W^{*}(R,D_{G}),E_{R}\right) =\left(
W^{*}(R,D_{R}),E_{D_{R}}^{D_{G}}\circ E\right) \otimes \left( D_{G},\mathbf{1%
}\right) $

\strut

$\ \ \ \ \ \ \ \ \ \ \ \ \ =\left( \underset{j=1}{\overset{N}{\,\,\,*_{D_{R}}%
}}\left( W^{*}(\{L_{w_{j}}\},D_{R}),E_{D_{R}}^{D_{G}}\circ E\right) \right)
\otimes \left( D_{G},\mathbf{1}\right) .$

\strut
\end{proof}

\strut

Now, we will provide the following fundamental examples ;

\strut \strut

\begin{example}
Let $G$ be a directed graph with $V(G)=\{v_{1},v_{2}\}$ and $%
E(G)=\{e=v_{1}ev_{2}\}.$ Let

\strut 

$\ \ \ \ \ \ \ \ \ \ \ \ \ \ \mathcal{F}=\{e\}$ \ \ and \ \ $%
R=\{L_{e},L_{e}^{*}\}.$

\strut 

We can construct the R-diagonal subalgebra $W^{*}(R,D_{G}).$ It is easy to
check that

\strut 

$\ \ \ \ \ \ 
\begin{array}{ll}
W^{*}\left( R,D_{G}\right)  & =W^{*}\left( \{L_{e},L_{e}^{*}\},D_{G}\right) 
\\ 
&  \\ 
& =W^{*}\left( \{L_{e},L_{e}^{*}\},D_{R}\right) \otimes D_{G} \\ 
&  \\ 
& =W^{*}\left( \{L_{e},L_{e}^{*},L_{v_{1}},L_{v_{2}}\}\right) .
\end{array}
$

\strut 

By the result from the final chapter, later, we can conclude that $%
W^{*}(R,D_{G})=W^{*}(G).$ So, the graph $W^{*}$-algebra $W^{*}(G)$ is same
as the R-diagonal subalgebra $W^{*}(R,D_{G})$ of it, where $R$ is consists
of all generators of $W^{*}(G)$ induced by the edge.
\end{example}

\strut

\begin{example}
Let $G$ be a directed graph with $V(G)=\{v_{1},v_{2}\}$ and $%
E(G)=\{e_{j}=v_{1}e_{j}v_{2}\}_{j=1}^{N}.$ Let $n\leq N.$ Then we have a
R-diagonal system,

\strut 

$\ \ \ \ \ \ \ R=\{L_{e_{1}},L_{e_{1}}^{*},...,L_{e_{n}},L_{e_{n}}^{*}\},$ $%
\forall n\leq N.$

\strut 

The R-diagonal subalgebra $W^{*}(R,D_{G})$ is trivially a $W^{*}$-subalgebra
of $W^{*}(G).$ Also,

\strut 

$\ \ \left( W^{*}(R,D_{G}),E_{R}\right) =\underset{j=1}{\overset{n}{%
\,\,*_{D_{G}}}}\left( W^{*}(\{L_{e_{j}}\},D_{G}),E_{j}\right) $

\strut 

$\ \ \ \ \ \ \ =\left( \underset{j=1}{\overset{n}{\,\,*_{D_{R}}}}\left(
W^{*}(\{L_{e_{j}}\},D_{R}),E_{j}\right) \right) \otimes \left( D_{G},\mathbf{%
1}\right) .$

\strut 
\end{example}

\strut \strut

\begin{example}
Let $G$ be a directed graph with $V(G)=\{v_{1},v_{2},v_{3}\}$ and

\strut 

$\ \ \ E(G)=\{e_{j}=v_{1}e_{j}v_{2},\,\,e_{k}^{\prime }=v_{2}e_{k}^{\prime
}v_{3}$

$\ \ \ \ \ \ \ \ \ \ \ \ \ \ \ \ \ \ \ \ \ \ \ \ \ \ \ \ \ \ \ \ \ \
:j=1,...,n$ \& $k=1,...,m\}.$

Take

$\ \ \ \ \ \ \ \ \ \ \ \ \ \ \ \ \ \ \ \ \ \mathcal{F}=\mathcal{F}_{1}%
\mathcal{\cup F}_{2},$

where

$\ \ \ \ \ \ \ \ \ \mathcal{F}_{1}=\{e_{1},...,e_{n}\}$ and $\mathcal{F}%
_{2}=\{e_{1}^{\prime },...,e_{m}^{\prime }\}.$

\strut 

If we construct the R-diagonal systems $R_{1}$ and $R_{2},$ induced by $%
\mathcal{F}_{1}$ and $\mathcal{F}_{2},$ respectively, then $R_{1}$ and $R_{2}
$ are free over $D_{G}$ in $\left( W^{*}(G),E\right) ,$ because they are
totally disjoint (See [16]). Therefore, the R-diagonal subalgebra generated
by $R_{1}\cup R_{2}$ is

\strut 

$\ W^{*}\left( R_{1}\cup R_{2},D_{G}\right)
=W^{*}(R_{1},D_{G})*_{D_{G}}W^{*}(R_{2},D_{G})$

\strut 

$\ \ \ =\left( \left( \underset{j=1}{\overset{n}{\,\,*_{D_{R_{1}}}}}\left(
W^{*}(\{L_{e_{j}}\},D_{R_{1}}),E_{e_{j}}\right) \right) \otimes \left( D_{G},%
\mathbf{1}\right) \right) $

\strut 

$\ \ \ \ \ \ \ *_{D_{G}}\left( \left( \underset{j=1}{\overset{m}{%
\,\,*_{D_{R_{2}}}}}\left( W^{*}(\{L_{e_{j}^{\prime
}}\},D_{R_{2}}),E_{e_{j}^{\prime }}\right) \right) \otimes \left( D_{G},%
\mathbf{1}\right) \right) .$
\end{example}

\strut

\strut

\strut

\strut

\section{Free Product Structures of $\left( W^{*}(G),E\right) $}

\strut

\strut

Throughout this chapter, let $G$ be a countable directed graph and $\left(
W^{*}(G),E\right) ,$ the graph $W^{*}$-probability space over the diagonal
subalgebra $D_{G}$. In this chapter, we will consider the building blocks of 
$W^{*}(G).$ Notice that if $w_{1}\neq w_{2}\in FP(G)$ and if $w_{1}w_{2}\in 
\Bbb{F}^{+}(G),$ then we can construct the $D_{G}$-valued random variable $%
L_{w_{1}w_{2}}$ which is same as $L_{w_{1}}L_{w_{2}}.$ Also notice that $%
L_{w_{1}}$ and $L_{w_{2}}$ are not free over $D_{G},$ in general. But under
the diagram-distinctness of $w_{1}$ and $w_{2},$ $L_{w_{1}}$ and $L_{w_{2}}$
are free over $D_{G}.$

\strut \strut \strut

\strut

\subsection{The $D_{G}$-Free Product Structures of $\left( W^{*}(G),E\right) 
$ I.}

\strut

\strut

In this section, we will consider the $D_{G}$-semicircular algebras and $%
D_{G}$-valued R-diagonal algebras in $W^{*}(G),$ more in detail. Recall that
the loop $l$ is basic if there is no other loop $w$ such that $l$ $=$ $%
w^{k}, $ for some $k$ $\in $ $\Bbb{N}\,\setminus \,\{1\}.$ Define

\strut

\begin{center}
$Loop(G)\overset{def}{=}\{l\in loop(G):l$ is basic$\}.$
\end{center}

\strut \strut \strut \strut

The following lemma is easily proved, by the very definition of basic loops ;

\strut

\begin{lemma}
Let $w\in loop(G)$ with $w=l^{k},$ for some $l\in Loop(G)$ and $k\in \Bbb{N}%
\,\setminus \,\{1\}.$ Then

\strut 

$\ \ \ \ \ \ \ \ \ \ W^{*}(\{L_{w}\},D_{G})\leq W^{*}(\{L_{l}\},D_{G}).$

$\square $
\end{lemma}

\strut

Remark that since $w$ and $l$ are not diagram-distinct, $L_{w}$ and $L_{l}$
are not free over $D_{G}$ in $\left( W^{*}(G),E\right) .$ In fact, $%
W^{*}(\{L_{l}\},D_{G})$ contains all $W^{*}(\{L_{w}\},D_{G}),$ if $l$ $\in $ 
$Loop(G)$ and $w$ $=$ $l^{k},$ $\forall $ $k$ $\in $ $\Bbb{N}.$

\strut \strut

\begin{proposition}
$\ \left( W^{*}\left( \{L_{l}:l\in loop(G)\},D_{G}\right) ,E\right) =%
\underset{l\in Loop(G)}{\,*_{D_{G}}}\left( W^{*}(\{L_{l}\},D_{G}),E\right) .$
\end{proposition}

\strut

\begin{proof}
Let

$\ \ \ \ \ \ \ \ \ \ \ \ \ \mathcal{L}=\{L_{l}:l\in loop(G)\},$

$\ \ \ \ \ \ \ \ \ \ \ \ \ \mathcal{L}_{l}=\{L_{l^{k}}:l\in Loop(G),\,k\in 
\Bbb{N}\}$

and

$\ \ \ \ \ \ \ \ \ \ \ \ \ \mathcal{L}_{0}=\{L_{l}:l\in Loop(G)\}.$

Then

$\ \ \ \ \ \ \ \ \mathcal{L}=\,\underset{l\in Loop(G)}{\cup }\mathcal{L}%
_{l}=\,\underset{L_{l}\in \mathcal{L}_{0}}{\cup }\,\left( \cup
_{k=1}^{\infty }\{L_{l}^{k}\}\right) .$

\strut

Thus

\strut

$\ \left( W^{*}\left( \{L_{l}:l\in loop(G)\},D_{G}\right) ,E\right) =\left(
W^{*}(\mathcal{L},D_{G}),E\right) $

\strut

$\ \ \ =\left( W^{*}(\,\underset{l\in Loop(G)}{\cup }\mathcal{L}%
_{l},D_{G}\,),E\right) $

\strut

$\ \ \ =\left( W^{*}(\,\underset{L_{l}\in \mathcal{L}_{0}}{\cup }\,\left(
\cup _{k=1}^{\infty }\{L_{l}^{k}\}\right) ,D_{G}\,),E\right) $

\strut

$\ \ \ =\underset{L_{l}\in \mathcal{L}_{0}}{\,\,*_{D_{g}}}\left(
W^{*}(\,\,\left( \cup _{k=1}^{\infty }\{L_{l}^{k}\}\right)
,D_{G}\,),E\right) $

\strut

by the fact that if $L_{l_{1}}\neq L_{l_{2}}$ in $\mathcal{L}_{0},$ then
they are free over $D_{G},$ by the diagram-distinctness of $l_{1}\neq
l_{2}\in Loop(G)$

\strut

$\ \ \ =\underset{L_{l}\in \mathcal{L}_{0}}{\,\,*_{D_{g}}}\left(
W^{*}(\,\,\left( \{L_{l}\}\right) ,D_{G}\,),E\right) ,$

\strut

since $W^{*}(\,\,\left( \cup _{k=1}^{\infty }\{L_{l}^{k}\}\right)
,D_{G})=W^{*}(\,\,\left( \{L_{l}\}\right) ,D_{G}\,).$ Therefore,

\strut

$\ \ \ \left( W^{*}\left( \{L_{l}:l\in loop(G)\},D_{G}\right) ,E\right) =%
\underset{l\in Loop(G)}{\,\,*_{D_{g}}}\left( W^{*}(\,\,\left(
\{L_{l}\}\right) ,D_{G}\,),E\right) .$

\strut
\end{proof}

\strut \strut

Finally, we can have the $D_{G}$-free product structure of the graph $W^{*}$%
-probability space $\left( W^{*}(G),E\right) $ over its diagonal subalgebra $%
D_{G}.$ By considering the $D_{G}$-freeness of generators of $W^{*}(G),$ we
can characterize the free product structure of $\left( W^{*}(G),E\right) .$

\strut

\begin{theorem}
Let $G$ be a countable directed graph and $\left( W^{*}(G),E\right) ,$ the
graph $W^{*}$-probability space over its diagonal subalgebra $D_{G}.$ Then

\strut 

$\ \ \ 
\begin{array}{ll}
\left( W^{*}(G),E\right)  & =\left( D_{G},E\right)  \\ 
&  \\ 
& \,\,\,\,\,\,\,*_{D_{G}}\left( \underset{l\in Loop(G)}{*_{D_{G}}}\left(
W^{*}(\{L_{l}\},D_{G}),E\right) \right)  \\ 
&  \\ 
& \,\,\,\,\,\,\,*_{D_{G}}\left( \underset{w\in loop^{c}(G)}{\,*_{D_{G}}}%
\left( W^{*}(\{L_{w}\},D_{G}),E\right) \right) .
\end{array}
$
\end{theorem}

\strut

\begin{proof}
Recall that $D_{G}$-valued random variables $L_{w_{1}}$ and $L_{w_{2}}$ are
free over $D_{G}$ in $\left( W^{*}(G),E\right) $ if and only if $w_{1}$ and $%
w_{2}$ are diagram-distinct. So, for any loop\thinspace $l$ and non-loop
finite path $w,$ $L_{l}$ and $L_{w}$ are free over $D_{G}.$ So,

\strut

$\ \ \ \ \ \ \ \ \ \ \ \ \ W^{*}\left( \{L_{l}:l\in loop(G)\},D_{G}\right) $

and

$\ \ \ \ \ \ \ \ \ \ \ W^{*}\left( \{L_{w}:w\in loop^{c}(G)\},D_{G}\right) $

\strut

are free over $D_{G}$ in $\left( W^{*}(G),E\right) .$ Denote the above
subalgebras by $\mathbf{L}$ and $\mathbf{R}$, respectively. Therefore, we
have that the free product space

\strut

$\ \ \ \ \ \ \ \ \ 
\begin{array}{ll}
D_{G} & *_{D_{G}}W^{*}\left( \{L_{l}:l\in loop(G)\},D_{G}\right) \\ 
&  \\ 
& *_{D_{G}}W^{*}\left( \{L_{w}:w\in loop^{c}(G)\},D_{G}\right)
\end{array}
\ $

\strut

is contained in $W^{*}(G).$ Since the generators of $W^{*}(G)$ and those of $%
D_{G}*_{D_{G}}\mathcal{L}$ $*_{D_{G}}\mathcal{R}$ are same, we can conclude
that

\strut

$\ \ \ \ \ \ \ \ \ \left( W^{*}(G),E\right) =\left( D_{G}*_{D_{G}}\mathbf{L}%
*_{D_{G}}\mathbf{R}\text{ },E\right) .$

\strut

But, by the previous proposition, we obtained that

\strut

$\ \ \ \ \ \ \ \ \ \left( \mathbf{L}\text{ },E\right) =\,\underset{l\in
Loop(G)}{\,*_{D_{G}}}\left( W^{*}(\{L_{l}\},D_{G}),E\right) .$

\strut

Now, we will observe that

\strut

$\ \ \ \ \ \ \ \ \ \left( \mathbf{R}\text{ },E\right) =\,\underset{w\in
loop^{c}(G)}{\,*_{D_{G}}}\left( W^{*}(\{L_{w}\},D_{G}),E\right) .$

\strut

Assume that $L_{w_{1}},L_{w_{2}}\in \left( \mathbf{R}\text{ },E\right) $ are
the generators. Then

\strut

$\ \ \ W^{*}\left( \{L_{w_{1}},L_{w_{2}}\},D_{G}\right)
=W^{*}(\{L_{w_{1}}\},D_{G})*_{D_{G}}W^{*}(\{L_{w_{2}}\},D_{G}).$

\strut

Therefore, we can get

\strut

$\ \ \ \ \ \ \ \ \ \left( \mathbf{R}\text{ },E\right) \leq \,\underset{w\in
loop^{c}(G)}{\,*_{D_{G}}}\left( W^{*}(\{L_{w}\},D_{G}),E\right) .$

\strut

The subalgebra inclusion ``$\geq $'' is clear. So,

\strut

$\ \ \ \ \left( W^{*}(G),E\right) =\left( D_{G}*_{D_{G}}\mathbf{L}*_{D_{G}}%
\mathbf{R}\text{ },E\right) $

\strut

$\ \ \ \ \ \ \ \ \ \ \ \ =\left( D_{G},E\right) *_{D_{G}}\left( \mathbf{L}%
,E\right) *_{D_{G}}\left( \mathbf{R}\text{ },E\right) $

\strut

$\ \ \ \ \ \ \ \ \ \ \ \ =\left( D_{G},E\right) *_{D_{G}}\left( \underset{%
l\in Loop(G)}{\,*_{D_{G}}}\left( W^{*}(\{L_{l}\},D_{G}),E\right) \right) $

\strut

$\ \ \ \ \ \ \ \ \ \ \ \ \ \ \ \ \ \ \ \ \ \ \ \ \ \ \ \ \ *_{D_{G}}\left( 
\underset{w\in loop^{c}(G)}{\,*_{D_{G}}}\left(
W^{*}(\{L_{w}\},D_{G}),E\right) \right) .$

\strut
\end{proof}

\strut \strut \strut \strut \strut

\strut

\strut \strut \strut

\subsection{\strut $D_{G}$-Free Building Blocks of $\left( W^{*}(G),E\right) 
$}

\strut

\strut

In this section, we will construct the $D_{G}$-free building blocks of the
graph $W^{*}$-probability space $\left( W^{*}(G),E\right) $ over its
diagonal subalgebra $D_{G}.$ Recall that

\strut

\begin{center}
$
\begin{array}{ll}
\left( W^{*}(G),E\right) = & \left( D_{G},E\right) \\ 
& \text{ \ }*_{D_{G}}\left( \underset{l\in Loop(G)}{\,*_{D_{G}}}\left(
W^{*}(\{L_{l}\},D_{G}),E\right) \right) \\ 
& \text{ \ }*_{D_{G}}\left( \underset{w\in loop^{c}(G)}{\,*_{D_{G}}}\left(
W^{*}(\{L_{w}\},D_{G}),E\right) \right) ,
\end{array}
$
\end{center}

$\strut \strut \ $

\strut where $Loop(G)$ is the collection of all basic loops contained in $%
loop(G).$ Notice that, even though we have a finite directed graph, $Loop(G)$
and $loop^{c}(G)$ may contain countably many elements. So, $\left(
W^{*}(G),E\right) $ is, in general, the $D_{G}$-free product of infinitely
many $D_{G}$-free $W^{*}$-subalgebras. But, in the final chapter, we will
show that this infinite free product of algebras can be contained in the
finite free product of algebras.

\strut \strut \strut \strut

\begin{definition}
Let $G$ be a countable directed graph and let $W^{*}(G)$ be the graph $W^{*}$%
-algebra. The diagonal subalgebra $D_{G}$ and $W^{*}$-subalgebras, $%
W^{*}(\{L_{l}\},$ $D_{G})$ and $W^{*}(\{L_{w}\},$ $D_{G})$ for all $l$ $\in $
$Loop(G)$ and $w$ $\in $ $loop^{c}(G)$ are $D_{G}$\textbf{-free building
blocks of }$W^{*}(G).$
\end{definition}

\strut \strut \strut \strut

As we observed in Chapter 1 and Chapter 2, we have that if $l$ is a loop,
then

\strut

(3.1) \ $\ \left( W^{*}(\{L_{l}\},D_{G}),E\right) =\left(
W^{*}(\{L_{l}\}),\,tr\right) \otimes \left( D_{G},\mathbf{1}\right) ,$

\strut

where $tr=E\mid _{W^{*}(\{L_{l}\})}$ is a tracial linear functional on $%
W^{*}(\{L_{l}\}).$ We also have that if $w$ is a non-loop finite path, then

\strut

(3.2)\ \ $\left( W^{*}(\{L_{w}\},D_{G}),E\right) =\left(
W^{*}(\{L_{w}\},D_{w}),E_{2}\right) \otimes \left( D_{G},\mathbf{1}\right) ,$

\strut

where $D_{w}=\overline{\Bbb{C}[\{L_{v_{1}},L_{v_{2}}:w=v_{1}wv_{2}\}]}^{w}$
is the $W^{*}$-subalgebra of $D_{G}$ and $E_{2}$ $=$ $E_{D_{w}}^{D_{G}}$ $%
\circ $ $E.$ By (3.2), we can get the following proposition which shows us
the vector space property of the non-loop $D_{G}$-free building blocks of $%
W^{*}(G).$

\strut

\begin{proposition}
Let $w\in loop^{c}(G)$ be a non-loop finite path and let $%
W^{*}(\{L_{w}\},D_{G})$ be the corresponding free building block. Then, as a
topological vector space,

\strut 

$\ \ W^{*}(\{L_{w}\},D_{G})=\overline{\Bbb{C}\{d,\,p(L_{w},L_{w}^{*}):d\in
D_{G},\,p\in \Bbb{C}_{1}[z_{1},z_{2}]\}}^{w},$

\strut 

where

\strut 

$\ \ \ \ \Bbb{C}_{1}[z_{1},z_{2}]\overset{def}{=}\{p\in \Bbb{C}%
[z_{1},z_{2}]:p(z_{1},z_{2})=\alpha _{0}+\alpha _{1}z_{1}+\alpha _{2}z_{2}\},
$

\strut 

for $\alpha _{0},\,\alpha _{1},\,\alpha _{2}\in \Bbb{C}.$ $\square $
\end{proposition}

\strut

\begin{proof}
In general, if $w\in FP(G)$ is a finite path, then the free building block $%
W^{*}(\{L_{w}\},D_{G})$ is a $W^{*}$-subalgebra of the graph $W^{*}$-algebra 
$W^{*}(G)$ such that

\strut

$\ \ \ W^{*}(\{L_{w}\},D_{G})=\overline{span\{d,\,p(L_{w},L_{w}^{*}):d\in
D_{G},\,p\in \Bbb{C}[z_{1},z_{2}]\}}^{w},$

\strut

as a topological vector space. Let $w\in loop^{c}(G)$ be a non-loop finite
path. Then $w^{k}\notin \Bbb{F}^{+}(G),$ for all $k\in \Bbb{N}\,\setminus
\,\{1\}.$ In other words, if $k\neq 1,$ then $w^{k}$ is not a admissible
finite path of the graph $G.$ Thus $L_{w}^{k}=L_{w}^{*\,\,k}=0_{D_{G}},$ for
all $k\in \Bbb{N}\,\setminus \,\{1\}.$ Therefore, if $q\in \Bbb{C}%
_{1}[z_{1},z_{2}],$ then

\strut

\ \ \ \ \ \ \ \ \ $q\left( L_{w},L_{w}^{*}\right) =\left\{ 
\begin{array}{lll}
q_{1}(L_{w},L_{w}^{*}) &  & \text{or} \\ 
\alpha \in \Bbb{C}\text{,} &  & 
\end{array}
\right. $

\strut

in general, where $q_{1}\in \Bbb{C}_{1}[z_{1},z_{2}].$ i.e.e,

\strut

(i) \ if $q(z_{1},z_{2})=\alpha _{0}+\sum_{k=2}^{\infty }\left( \alpha
_{k}^{1}z_{1}^{k}+\alpha _{k}^{2}z_{2}^{k}\right) ,$ then

\strut

$\ \ \ \ \ \ \ \ \ \ \ \ \ \ \ \ \ q\left( L_{w},L_{w}^{*}\right) =\alpha
_{0}.$

\strut

(ii) if $q(z_{1},z_{2})=\alpha _{0}+\sum_{k=1}^{\infty }\left( \alpha
_{k}^{1}z_{1}^{k}+\alpha _{k}^{2}z_{2}^{k}\right) ,$ then

\strut

$\ \ \ \ \ \ \ \ \ q\left( L_{w},L_{w}^{*}\right) =\alpha _{0}+\left( \alpha
_{1}^{1}L_{w}+\alpha _{1}^{2}L_{w}^{*}\right) .$

\strut

So, if we define $q_{1}\in \Bbb{C}_{1}[z_{1},z_{2}]$ by

\strut

$\ \ \ \ \ \ \ \ \ \ \ q_{1}(z_{1},z_{2})=\alpha _{0}+\left( \alpha
_{1}^{1}z_{1}+\alpha _{1}^{2}z_{2}\right) ,$

then

$\ \ \ \ \ \ \ \ \ \ \ \ \ \ \ q\left( L_{w},L_{w}^{*}\right) $ $=$ $%
q_{1}\left( L_{w},L_{w}^{*}\right) .$

\strut
\end{proof}

\strut

\strut

\subsection{The $D_{G}$-Free Product Structure of $\left( W^{*}(G),E\right) $
II.}

\strut

\strut

By Section 3.1 and by (3.1) and (3.2), we can get the following theorem;

\strut

\begin{theorem}
Let $G$ be a countable directed graph and $\left( W^{*}(G),E\right) ,$ the
graph $W^{*}$-probability space over its diagonal subalgebra $D_{G}.$ Then

\strut 

$\ 
\begin{array}{ll}
\left( W^{*}(G),E\right)  & =D_{G} \\ 
&  \\ 
& \,\,\,*_{D_{G}}\left( \underset{l\in Loop(G)}{*_{D_{G}}}\left( \left(
W^{*}(\{L_{l}\}),tr_{l}\right) \otimes \left( D_{G},\mathbf{1}\right)
\right) \right)  \\ 
&  \\ 
& \,\,\,*_{D_{G}}\left( \underset{w\in loop^{c}(G)}{*_{D_{G}}}\left( \left(
W^{*}(\{L_{w}\},D_{w}),E_{w}\right) \otimes \left( D_{G},\mathbf{1}\right)
\right) \right) ,
\end{array}
$

\strut 

where $tr_{l}=E_{D_{l}}^{D_{G}}\circ E$ is a trace on $W^{*}(\{L_{l}\})$ and 
$E_{w}=E_{D_{w}}^{D_{G}}\circ E$ is a conditional expectation from $W^{*}(G)$
onto $D_{w}$ ($E_{B}^{A}$ means the conditional expectation from $A$ onto $B$%
) and

\strut 

$\ \ \ \ \ \ \ \ \ \ \ \ \ \ \ \ \ D_{l}$ $=$ $\overline{\Bbb{C}%
[\{L_{v}:l=vlv\}]}^{w}=\Bbb{C}$

and

$\ \ \ \ \ \ \ D_{w}=\overline{\Bbb{C}[\{L_{v_{1}},L_{v_{2}}:w=v_{1}wv_{2}\}]%
}^{w}=\Delta _{2}.$

$\square $
\end{theorem}

\strut

\strut \strut

\strut \strut \strut

\section{More About the Free Product Structure of $\left( W^{*}(G),E\right) $%
}

\strut \strut

\strut

In this chapter, we will complete to observe the free product structure of
the graph $W^{*}$-probability spaces. Let $G$ be a countable directed graph
and let $\left( W^{*}(G),E\right) $ be the graph $W^{*}$-probability space
over its diagonal subalgebra $D_{G}.$ In Chapter 3, we showed that

\strut

\begin{center}
$
\begin{array}{ll}
\left( W^{*}(G),E\right) = & \left( D_{G},\mathbf{1}\right) \\ 
&  \\ 
& \,\,*_{D_{G}}\left( \underset{l\in Loop(G)}{*_{D_{G}}}\left( \left(
W^{*}(\{L_{l}\}),tr\right) \otimes \left( D_{G},\mathbf{1}\right) \right)
\right) \\ 
&  \\ 
& \,\,*_{D_{G}D_{G}}\left( \underset{w\in loop^{c}(G)}{*_{D_{G}}}\left(
\left( W^{*}(\{L_{w}\},D_{w}),E_{w}\right) \otimes \left( D_{G},\mathbf{1}%
\right) \right) \right) .
\end{array}
$
\end{center}

\strut

In the previous chapter, we emphasize the roles of free building blocks and
tried to consider each free building block. By using the characterization of
the free building blocks, we could get the above free product structure of
the graph $W^{*}$-probability space $\left( W^{*}(G),E\right) .$ Without
considering the structure of each free building blocks of $\left(
W^{*}(G),E\right) ,$ by Section 3.1, we can rewrite the above formula as

\strut

(4.1)

\begin{center}
$
\begin{array}{ll}
\left( W^{*}(G),E\right) & =(D_{G},E) \\ 
&  \\ 
& \,\,\,\,*_{D_{G}}\left( \underset{l\in Loop(G)}{*_{D_{G}}}\left(
W^{*}(\{L_{l}\},D_{G}),E\right) \right) \\ 
&  \\ 
& \,\,\,\,*_{D_{G}}\left( \underset{w\in loop^{c}(G)}{*_{D_{G}}}\left(
W^{*}(\{L_{w}\},D_{G}),E\right) \right) .
\end{array}
$
\end{center}

\strut

In this chapter, we will show that

\strut (4.2)

\begin{center}
$
\begin{array}{ll}
\left( W^{*}(G),E\right) & =(D_{G},E) \\ 
&  \\ 
& \,\,\,\,*_{D_{G}}\left( \underset{l\in ELoop(G)}{*_{D_{G}}}\left(
W^{*}(\{L_{l}\},D_{G}),E\right) \right) \\ 
&  \\ 
& \,\,\,\,*_{D_{G}}\left( \underset{w\in Eloop^{c}(G)}{*_{D_{G}}}\left(
W^{*}(\{L_{w}\},D_{G}),E\right) \right) .
\end{array}
$
\end{center}

\strut

where\strut

\begin{center}
$ELoop(G)\overset{def}{=}E(G)\cap Loop(G)$
\end{center}

and

\begin{center}
$Eloop^{c}(G)\overset{def}{=}E(G)\cap loop^{c}(G).$
\end{center}

\strut

Equivalently, we will show that

\strut

(4.3) $\ \left( W^{*}(G),E\right) =\left( D_{G},\mathbf{1}\right)
*_{D_{G}}\left( \underset{e\in E(G)}{\,*_{D_{G}}}\left(
W^{*}(\{L_{e}\},D_{G}),E\right) \right) .$

\strut

First, we will concentrate on proving the formula (4.1) is equivalent to the
formula (4.2). For the convenience, define

\strut

\begin{center}
$\mathcal{L}_{edge}(G)\overset{def}{=}\underset{l\in ELoop(G)}{*_{D_{G}}}%
W^{*}(\{L_{l}\},D_{G}),$
\end{center}

\strut

\begin{center}
$\mathcal{L}_{edge}^{c}(G)\overset{def}{=}\underset{w\in Eloop^{c}(G)}{%
*_{D_{G}}}W^{*}(\{L_{w}\},D_{G}),$
\end{center}

\strut

\begin{center}
$\mathcal{L}(G)\overset{def}{=}\underset{l\in Loop(G)}{*_{D_{G}}}%
W^{*}(\{L_{l}\},D_{G}),$
\end{center}

and

\begin{center}
$\mathcal{L}^{c}(G)\overset{def}{=}\underset{w\in loop^{c}(G)}{*_{D_{G}}}%
W^{*}(\{L_{w}\},D_{G}).$
\end{center}

\strut

\begin{theorem}
Let $G$ be a countable directed graph and let $\left( W^{*}(G),E\right) $ be
the graph $W^{*}$-probability space over its diagonal subalgebra $D_{G}.$
Then

\strut \strut 

$\ \ \ \left( W^{*}(G),E\right) =\left( D_{G},E\right) *_{D_{G}}\left( 
\mathcal{L}_{edge}(G),E\right) *_{D_{G}}\left( \mathcal{L}%
_{edge}^{c}(G),E\right) .$
\end{theorem}

\strut

\begin{proof}
Let $ELoop(G)$ $=$ $E(G)$ $\cap $ $Loop(G)$ and $Eloop^{c}(G)$ $=$ $E(G)$ $%
\cap $ $loop^{c}(G).$ By (4.1), \strut

\strut

$\ \ \ \ \ \ \ \left( W^{*}(G),E\right) =(D_{G},E)*_{D_{G}}\left( \mathcal{L}%
(G),E\right) *_{D_{G}}\left( \mathcal{L}^{c}(G),E\right) .$

\strut

Since $Eloop(G)\subset Loop(G)$ and $Eloop^{c}(G)\subset loop^{c}(G),$ we
have that

\strut

$\ \ \ \ \ \ \ \ \ \ \ \mathcal{L}_{edge}(G)\leq \mathcal{L}(G)$ \ \ and \ \ 
$\mathcal{L}_{edge}^{c}(G)\leq \mathcal{L}^{c}(G).$

\strut

Therefore, we have the following subalgebra inclusion ``$\geq $'' ;

\strut

(4.4)$\ \left( W^{*}(G),E\right) \geq \left( D_{G},E\right) *_{D_{G}}\left( 
\mathcal{L}_{edge}(G),E\right) *_{D_{G}}\left( \mathcal{L}%
_{edge}^{c}(G),E\right) .$

\strut

So, it suffices to show that we have the reverse subalgebra inclusion ``$%
\leq $''.

\strut

(\textbf{Case I}) Assume that $l\in Loop(G).$ If $l\in ELoop(G),$ then

\strut

(4.5)$\ \ \ \ \ \ \ \ \ \ \ \ \ \ \ W^{*}(\{L_{l}\},D_{G})<\mathcal{L}%
_{edge}(G).$

\strut

If $l=e_{1},...,e_{n}$ with $e_{1},...,e_{n}\in E(G),$ such that $e_{j}\in
Eloop^{c}(G),$ for all $j=1,...,n,$ $n>1,$ then we have that

\strut

$\ \ \ \ \ \ \ \ \ \ \ \ \ W^{*}(\{L_{l}\},D_{G})\leq \underset{j=1}{%
\overset{n}{\,*_{D_{G}}}}W^{*}(\{L_{e_{j}}\},D_{G}).$

\strut

Therefore, by the assumption that $e_{j}\in Eloop^{c}(G),$ we have that

\strut

(4.6) \ \ \ \ \ \ \ \ \ \ \ \ \ \ \ $W^{*}\left( \{L_{l}\},D_{G}\right) <%
\mathcal{L}_{edge}^{c}(G).$

\strut

By (4.5) and (4.6), we can conclude that if $l$ $\in $ $Loop(G),$ then the
subalgebra $W^{*}(\{L_{l}\},$ $D_{G})$ of $W^{*}(G)$ is the subalgebra of $%
\mathcal{L}_{edge}(G)$ $*_{D_{G}}$ $\mathcal{L}_{edge}^{c}(G).$ i.e.,

\strut

(4.7) \ $\left( W^{*}(\{L_{l}\},D_{G})\right) \leq \left( \mathcal{L}%
_{edge}(G),E\right) *_{D_{G}}\left( \mathcal{L}_{edge}^{c}(G),E\right) ,$

\strut

for all $l\in Loop(G).$ Therefore,

\strut

(4.8) \ \ $\ \ \ \left( \mathcal{L}(G),E\right) \leq \left( \mathcal{L}%
_{edge}(G),E\right) *_{D_{G}}\left( \mathcal{L}_{edge}^{c}(G),D_{G}\right) .$

\strut

(\textbf{Case II}) Now, assume that $w\in loop^{c}(G).$ Suppose that $w$ $%
\in $ $Eloop^{c}(G).$ Then, clearly,

\strut

(4.9) \ \ \ \ \ \ \ \ \ \ $W^{*}(\{L_{w}\},D_{G})<\mathcal{L}_{edge}^{c}(G).$

\strut

Now, assume that $w=e_{1}...e_{k}\in loop^{c}(G)$ with $e_{1},...,e_{k}\in
E(G)$ are edges satisfying that the initial vertex of $e_{1}$ and the final
vertex of $e_{k}$ are different. Then

\strut

(4.10)\ \ \ \ \ \ $W^{*}(\{L_{w}\},D_{G})\leq \underset{j=1}{\overset{k}{%
\,\,\,*_{D_{G}}}}W^{*}(\{L_{e_{j}}\},D_{G}).$

\strut

This also shows that

\strut

(4.11) \ \ \ \ \ \ \ \ \ \ $W^{*}\left( \{L_{w}\},D_{G}\right) <\mathcal{L}%
_{edge}^{c}(G).$

\strut

By (4.10) and (4.11), we can conclude that if $w$ $\in $ $loop^{c}(G),$ then
the subalgebra $W^{*}(\{L_{w}\},D_{G})$ of $W^{*}(G)$ is the subalgebra of $%
\mathcal{L}_{edge}^{c}(G).$ i.e.,

\strut

(4.12) \ \ \ \ \ \ $\left( W^{*}(\{L_{w}\},D_{G}),E\right) \leq \left( 
\mathcal{L}_{edge}^{c}(G),E\right) ,$

\strut

for all $w\in loop^{c}(G).$ Therefore, we have that

\strut

(4.13) \ \ \ \ \ \ \ \ \ $\left( \mathcal{L}^{c}(G),E\right) \leq \left( 
\mathcal{L}_{edge}^{c}(G),E\right) .$

\strut

As we considered in the previous two cases, we can get that

\strut

(4.14) \ \ \ \ $\mathcal{L}(G)*_{D_{G}}\mathcal{L}^{c}(G)\leq \mathcal{L}%
_{edge}(G)*_{D_{G}}\mathcal{L}_{edge}^{c}(G).$

\strut

Therefore, by the relation (4.14), we can conclude that

\strut

$\ \ \ \left( W^{*}(G),E\right) \leq \left( D_{G},E\right) *_{D_{G}}\left( 
\mathcal{L}_{edge}(G),E\right) *_{D_{G}}\left( \mathcal{L}%
_{edge}^{c}(G),E\right) .$

\strut
\end{proof}

\strut

The above theorem provides us that

\strut

\begin{center}
$
\begin{array}{ll}
\left( W^{*}(G),E\right) & =(D_{G},E) \\ 
&  \\ 
& \,\,\,\,*_{D_{G}}\left( \underset{l\in E(G)\cap Loop(G)}{*_{D_{G}}}\left(
W^{*}(\{L_{l}\},D_{G}),E\right) \right) \\ 
&  \\ 
& \,\,\,\,*_{D_{G}}\left( \underset{w\in E(G)\cap loop^{c}(G)}{*_{D_{G}}}%
\left( W^{*}(\{L_{w}\},D_{G}),E\right) \right) .
\end{array}
$
\end{center}

\strut

Therefore, we can get the following simple free product structure of a graph 
$W^{*}$-probability space $\left( W^{*}(G),E\right) $ ;

\strut

\begin{corollary}
Let $G$ be a countable directed graph and let $\left( W^{*}(G),E\right) $ be
the graph $W^{*}$-probability space over its diagonal subalgebra $D_{G}.$
Then $\left( W^{*}(G),E\right) $ has the following free product structure ;

\strut 

$\ \ \ \left( W^{*}(G),E\right) =\left( D_{G},E\right) *_{D_{G}}\left( 
\underset{e\in E(G)}{\,*_{D_{G}}}\left( W^{*}(\{L_{e}\},D_{G}),E\right)
\right) .$
\end{corollary}

\strut

\begin{proof}
Notice that the edge set $E(G)$ of the graph $G$ satisfies that

\strut

$\ \ \ \ \ \ \ E(G)=\left( E(G)\cap Loop(G)\right) \cup \left(
Eloop^{c}(G)\right) .$

\strut

Assume that if $e\in E(G)$ is a loop-edge in $loop(G)$, then $e$ is a basic
loop in $Loop(G).$ i.e.e,

\strut

$\ \ \ \ \ \ \ \ \ \ \ E(G)\cap Loop(G)=Eloop(G).$

\strut

Thus we have that

\strut

$\ \ \ 
\begin{array}{ll}
E(G) & =\left( E(G)\cap Loop(G)\right) \cup \left( E(G)\cap
loop^{c}(G)\right) \\ 
&  \\ 
& =\left( E(G)\cap loop(G)\right) \cup \left( E(G)\cap loop^{c}(G)\right) \\ 
&  \\ 
& =E(G)\cap \left( loop(G)\cup loop^{c}(G)\right) \\ 
&  \\ 
& =E(G)\cap FP(G)=E(G).
\end{array}
$

\strut \strut \strut
\end{proof}

\strut

The above theorem and corollary shows us that the free product structure of
a graph $W^{*}$-probability space $\left( W^{*}(G),E\right) $ is totally
depending on the admissibility on edges of the graph $G.$

\strut \strut \strut

In the rest of this section, we will consider the several examples ;

\strut \strut \strut \strut

\begin{example}
Let $G$ be a directed one-vertex graph with $V(G)=\{v\}$ and $%
E(G)=\{l_{1},...,l_{N}\},$ where $l_{j}$'s are all loop edges, \ $j=1,...,N.$
Then, by the previous theorem, we can have that

\strut 

$\ \ \ \left( W^{*}(G),E\right) =\left( D_{G},E\right) *_{D_{G}}\left( 
\underset{j=1}{\overset{N}{\,*_{D_{G}}}}\left(
W^{*}(\{L_{e_{j}}\},D_{G}),E\right) \right) .$

\strut 

Notice that $D_{G}=\Bbb{C}$ and $E=tr,$ where $tr$ is the canonical trace,
induced by the given conditional expectation. Therefore,

\strut 

$\ \ \ \ \ \ \ \ \left( W^{*}(G),tr\right) =\underset{j=1}{\overset{N}{%
\,*_{D_{G}}}}\left( W^{*}(\{L_{e_{j}}\}),tr\right) .$

\strut 

Trivially, it contains the free group factor $\left( L(F_{N}),\tau \right) ,$
where $\tau =tr\mid _{A},$ where

\strut 

$\ \ \ \ \ \ \ \ \ \ \ A=W^{*}\left(
\{L_{e_{j}}+L_{e_{j}}^{*}:j=1,...,N\}\right) .$
\end{example}

\strut \strut

\begin{example}
Suppose that we have a directed graph $G$ with

\strut 

$\ \ \ \ \ \ \ \ \ V(G)=\{v_{1},v_{2}\}$ \ and $E(G)=%
\{l_{1}^{1},l_{2}^{1},l_{1}^{2},l_{2}^{2},l_{3}^{2},e\},$

\strut where

\ \ $\ \ \ \ \ \ \ \ \ \ l_{j}^{1}=v_{1}l_{j}^{1}v_{1},$ for \ $j=1,2,$ \ \
\ $e=v_{1}ev_{2}$

and

$\ \ \ \ \ \ \ \ \ \ \ \ \ \ \ l_{j}^{2}=v_{2}l_{j}^{2}v_{2},$ for \ $%
j=1,2,3.$

\strut 

Then $D_{G}=\Delta _{2}\subset M_{2}(\Bbb{C})$ and

\strut 

$\ \ \ \ \ 
\begin{array}{ll}
\left( W^{*}(G),E\right)  & =D_{G}*_{D_{G}}\left( \underset{i=1}{\overset{2}{%
*_{D_{G}}}}\left( W^{*}(\{L_{l_{i}^{1}}\})\otimes D_{G},tr\otimes \mathbf{1}%
\right) \right)  \\ 
&  \\ 
& \,\,\,\,\,\,\,\,\,\,\,\,*_{D_{G}}\left( \underset{j=1}{\overset{3}{%
*_{D_{G}}}}\left( W^{*}(\{L_{l_{j}^{2}}\})\otimes D_{G},tr\otimes \mathbf{1}%
\right) \right)  \\ 
&  \\ 
& \,\,\,\,\,\,\,\,\,\,\,\,*_{D_{G}}\left( \left( W^{*}(\{L_{e}\},\Delta
_{2}),E_{2}\right) \otimes \left( D_{G},\mathbf{1}\right) \right) .
\end{array}
$

\strut \strut \strut 

This shows that our graph $W^{*}$-probability space $\left(
W^{*}(G),E\right) $ contains free group factors

\strut 

$\ \ \ \ \ \ \ \ \ \ \ \ \ L(F_{2})\leq *_{i=1}^{2}\left(
W^{*}(\{L_{l_{i}^{1}}\}),tr\right) $

and

$\ \ \ \ \ \ \ \ \ \ \ \ \ L(F_{3})\leq *_{j=1}^{3}\left(
W^{*}(\{L_{l_{j}^{2}}\}),tr\right) $

\strut 

(See Section 4.1). Notice that these free group factors $L(F_{2})$ and $%
L(F_{3})$ are free \textbf{over }$D_{G}$ in $\left( W^{*}(G),E\right) .$
Therefore, $\left( W^{*}(G),E\right) $ contains $L(F_{2})*_{D_{G}}L(F_{3}).$
Remark that

\strut 

$\ \ \ \ \ \ \ L(F_{2})*_{D_{G}}L(F_{3})\neq L(F_{2})*L(F_{3})=L(F_{5}).$
\end{example}

\strut

\begin{example}
Let $C_{N}$ be the circular graph with $V(C_{N})=\{v_{1},...,v_{N}\}$ and

\strut 

$\ \ \ E(C_{N})=\{e_{j}:e_{j}=v_{j}e_{j}v_{j+1},$ $j=1,...,N-1,$ $%
e_{N}=v_{N}e_{N}v_{1}\}.$

\strut 

Then $D_{G}=\Delta _{N}.$ Thus we have that

\strut 

$\ \ \ \ \ \ \ \ \ \ \ \ \ \ \ \ \ \ \ \ \ \ \ \ \
Loop(G)=\{l:=e_{1}...e_{N}\}$

equivalently,

$\ \ \ \ \ \ \ \ \ \ \ \ \ \ \ \ \ \ \ loop^{c}(G)=\{w\in FP(G):w\neq l^{k},$
$k\in \Bbb{N}\}.$

\strut 

So, for the canonical conditional expectation $E,$ we have that

\strut 

$\ \ \left( W^{*}(C_{N}),E\right) =\left( \Delta _{N},\mathbf{1}\right) $

\strut 

$\ \ \ \ \ \ \ \ \ *_{\Delta _{N}}\left( W^{*}(\{L_{l}\})\otimes \Delta
_{N},E\otimes \mathbf{1}\right) $

\strut 

$\ \ \ \ \ \ \ \ \ *_{\Delta _{N}}\left( \underset{w\in loop^{c}(G)}{%
\,*_{D_{G}}}\left( W^{*}(\{L_{w}\},\Delta _{2})\otimes \Delta
_{N},E_{2}\otimes \mathbf{1}\right) \right) $

\strut 

$\ \ \ =\left( \Delta _{N},E\right) *_{\Delta _{N}}\left( \underset{j=1}{%
\overset{N}{\,*_{\Delta _{N}}}}\left( W^{*}(\{L_{e_{j}}\},\Delta
_{N}),E\right) \right) .$
\end{example}

\strut \bigskip \strut \strut \strut

\strut \strut

\strut

\begin{quote}
\textbf{Reference}

\strut

\strut

{\small [1] \ \ A. Nica, R-transform in Free Probability, IHP course note,
available at www.math.uwaterloo.ca/\symbol{126}anica.}

{\small [2]\strut \ \ \ A. Nica and R. Speicher, R-diagonal Pair-A Common
Approach to Haar Unitaries and Circular Elements, (1995), www
.mast.queensu.ca/\symbol{126}speicher.\strut }

{\small [3] \ }$\ ${\small B. Solel, You can see the arrows in a Quiver
Operator Algebras, (2000), preprint}

{\small \strut [4] \ \ A. Nica, D. Shlyakhtenko and R. Speicher, R-cyclic
Families of Matrices in Free Probability, J. of Funct Anal, 188 (2002),
227-271.}

{\small [5] \ \ D. Shlyakhtenko, Some Applications of Freeness with
Amalgamation, J. Reine Angew. Math, 500 (1998), 191-212.\strut }

{\small [6] \ \ D.Voiculescu, K. Dykemma and A. Nica, Free Random Variables,
CRM Monograph Series Vol 1 (1992).\strut }

{\small [7] \ \ D. Voiculescu, Operations on Certain Non-commuting
Operator-Valued Random Variables, Ast\'{e}risque, 232 (1995), 243-275.\strut 
}

{\small [10]\ D. Shlyakhtenko, A-Valued Semicircular Systems, J. of Funct
Anal, 166 (1999), 1-47.\strut }

{\small [10]\ D.W. Kribs and M.T. Jury, Ideal Structure in Free Semigroupoid
Algebras from Directed Graphs, preprint}

{\small [10]\ D.W. Kribs and S.C. Power, Free Semigroupoid Algebras, preprint%
}

{\small [11]\ I. Cho, Amalgamated Boxed Convolution and Amalgamated
R-transform Theory, (2002), preprint.}

{\small [12] I. Cho, The Tower of Amalgamated Noncommutative Probability
Spaces, (2002), Preprint.}

{\small [13] I. Cho, Free Perturbed R-transform Theory, (2003), Preprint.}

{\small [14]\ I. Cho, Compatibility of a Noncommutative Probability Space
and a Noncommutative Probability Space with Amalgamation, (2003), Preprint}

{\small [15] I. Cho, An Example of Scalar-Valued Moments, Under
Compatibility, (2003), Preprint.}

{\small [16] I. Cho, Graph }$W^{\ast }$-{\small Probability Theory, (2004),
Preprint.}

{\small [17] I. Cho, Random Variables in Graph }$W^{*}$-{\small Probability
Spaces, (2004), Preprint. }

{\small [18] I. Cho, Amalgamated Semicircular Systems in Graph }$W^{*}$%
{\small -Probability Spaces, (2004), Preprint. }

{\small [19] I. Cho, Amalgamated R-diagonal Pairs, (2004), Preprint. }

{\small [20] I. Cho, Compressed Random Variables in Graph }$W^{*}${\small %
-Probability Spaces, (2004), Preprint.}

{\small [21] P.\'{S}niady and R.Speicher, Continous Family of Invariant
Subspaces for R-diagonal Operators, Invent Math, 146, (2001) 329-363.}

{\small [22] R. Speicher, Combinatorial Theory of the Free Product with
Amalgamation and Operator-Valued Free Probability Theory, AMS Mem, Vol 132 ,
Num 627 , (1998).}

{\small [23] R. Speicher, Combinatorics of Free Probability Theory IHP
course note, available at www.mast.queensu.ca/\symbol{126}speicher.\strut }
\end{quote}

\end{document}